\documentclass[12pt]{article}
\usepackage{amssymb,amsthm}

\newtheorem{definition}{Definition}[section]
\newtheorem{theorem}[definition]{Theorem}
\newtheorem{corollary}[definition]{Corollary}
\newtheorem{lemma}[definition]{Lemma}

\newtheorem{example}[definition]{Example}

\def\R{{\mathbb R}}
\def\C{{\mathbb C}}

\def\Mat{\rm Mat}

\begin{document}

\title{The Pseudo Cosine Sequences  \\
of a Distance-Regular Graph}
\author{Arlene A. Pascasio \qquad\qquad\qquad Paul M. Terwilliger}
\date{July 11,2003}
\maketitle

\begin{abstract}
Let $\Gamma$ denote a distance-regular graph with diameter $D\ge 3$, valency $k$, and intersection
numbers $a_i, b_i, c_i.$  By a \emph{pseudo cosine sequence} 
of $\Gamma$ we mean
a sequence of real numbers $\sigma_0,\sigma_1, \dots, \sigma_D$
such that $\sigma_0=1$ and
$
c_i\sigma_{i-1}+a_i\sigma_i+b_i\sigma_{i+1}=k\sigma_1\sigma_i$ for $ 0 \le i
\le D-1.$
Let
$\sigma_0, \sigma_1, \ldots, \sigma_D$  and
$\rho_0, \rho_1, \ldots, \rho_D$ denote pseudo cosine sequences
of $\Gamma$. We  say this pair of
sequences is \emph{tight} whenever
 $\sigma_0\rho_0, \sigma_1\rho_1, \ldots, \sigma_D\rho_D$ is a
pseudo cosine sequence of $\Gamma$. In this paper, we determine all the tight pairs
of pseudo cosine sequences of $\Gamma.$

\end{abstract}

{\small
\noindent AMS Classification: 05E50 \\
\noindent Keywords: distance-regular graph, association scheme, cosine sequence, pseudo cosine sequence.
}

\section{Introduction}

Let $\Gamma=(X,R)$ denote a distance-regular graph with
diameter $D\geq 3$, valency $k$, and intersection numbers 
$a_i, b_i, c_i$
(see Section 2 for formal definitions).
Let $\theta $ and $\theta'$ denote 
eigenvalues of $\Gamma$ other than $k$. In  
\cite{jkt}, Juri\v si\'c, Koolen and Terwilliger showed
\begin{equation}
\Biggl(\theta + \frac{k}{a_1+1}\Biggr)
\Biggl(\theta' + \frac{k}{a_1+1}\Biggr) \geq -\frac{ka_1b_1}{(a_1+1)^2}.
\label{eq:ththpintro} 
\end{equation}
The case of equality in
(\ref{eq:ththpintro}) has received a lot of attention.
This case has been characterized using
 the intersection
numbers \cite[Theorem 10.2]{jkt}, the 1-homogeneous
property \cite[Theorem 11.7(i)]{jkt}, and the local structure 
\cite[Theorem 12.6]{jkt}. 
See also \cite{go2, jk5, jk, jk0, jk2, jk3, jurter, maclean1, maclean2,
maclean3, aap1, aap2, aap3, aap4, tomiyama, terweng}.
In \cite{aap1}, Pascasio characterized  equality in
(\ref{eq:ththpintro}) 
using cosine sequences. The idea is as follows.
Let $\theta$ denote an eigenvalue of $\Gamma$.
Recall the {\it cosine sequence} for $\theta$ is 
the sequence of real numbers 
$\sigma_0, \sigma_1, \ldots, \sigma_D$ 
such that $\sigma_0=1$ and $c_i\sigma_{i-1}+a_i\sigma_i+b_i\sigma_{i+1}
=\theta \sigma_i$ for $0 \leq i \leq D$, where $\sigma_{-1}$ and
$\sigma_{D+1}$ 
are indeterminate. Let 
$\theta $ and $\theta'$ denote eigenvalues of $\Gamma $ other
than $k$.
Let $\sigma_0, \sigma_1, \ldots, \sigma_D$ 
(resp. $\rho_0, \rho_1, \ldots, \rho_D$) denote the cosine
sequence for $\theta$ 
(resp. $\theta'$). 
Then by \cite[Theorem 3.4]{aap1},
 $\theta, \theta'$ satisfy
(\ref{eq:ththpintro}) with equality if and only if 
$\sigma_0\rho_0, \sigma_1\rho_1, \ldots, \sigma_D\rho_D$ is a cosine
sequence.
 
\medskip
The concept of a {\it pseudo cosine sequence}
was recently introduced
 \cite{jurter, terweng}.
%
This concept is defined as follows. Let $\theta$ denote a real number.
By the {\it pseudo cosine sequence} for $\theta$ we mean
the sequence of real numbers  
 $\sigma_0, \sigma_1, \ldots, \sigma_D$ 
such that $\sigma_0=1$ and $c_i\sigma_{i-1}+a_i\sigma_i+b_i\sigma_{i+1}
=\theta \sigma_i$ for $0 \leq i \leq D-1$, where $\sigma_{-1}$ is
indeterminate.
Pseudo cosine sequences have been used to describe
certain modules for the subconstituent algebra \cite{terweng}.
They also arise in connection with the
 pseudo 1-homogeneous property \cite{jurter}. 
Given our comments in the previous paragraph it is natural
to consider the following situation.
Let 
$\sigma_0, \sigma_1, \ldots, \sigma_D$  and
 $\rho_0, \rho_1, \ldots, \rho_D$
denote pseudo cosine sequences of $\Gamma$. We say this pair
of sequences is {\it tight} whenever
$\sigma_0\rho_0, \sigma_1\rho_1, \ldots, \sigma_D\rho_D$
is a pseudo cosine sequence of $\Gamma$.

\medskip
In this paper we find all the tight pairs of
pseudo cosine sequences of $\Gamma$.
We break the argument into the following cases: 
(i) $a_i= 0$ for $0 \leq i \leq D-1$; (ii)
$a_i= 0$ for $0 \leq i \leq D-2$ and $a_{D-1}\not=0$;
(iii) $a_1=0$ and there exists $i$ $(2 \leq i \leq D-2)$ such
that $a_i\not=0$;
(iv) $a_1 \not=0$.
Our results for these cases are given in
Theorem
\ref{thm:bpclass},
Theorem
\ref{thm:case2},
Theorem 
\ref{thm:case3},
and Theorem
\ref{th:mainth} respectively.

\medskip
Our methods in this paper are purely algebraic.
Consequently our main results 
probably hold for
$P$-polynomial character algebras. We will pursue this in a future
paper.

\section{Preliminaries}

In this section, we review some definitions and basic concepts.
For more
background information, the reader may refer to the books of Bannai and
Ito \cite {bi}, or Brouwer, Cohen and Neumaier \cite {bcn}.\\

Let $\Gamma=(X,R)$ denote a finite, undirected, connected graph
without loops or multiple edges, with vertex set $X$, edge set $R$, path-length distance
function $\partial$ and diameter $D:={\rm max}\ \{\partial(x,y)| x,y \in X\}$.
Let $k$ denote a nonnegative integer. We say $\Gamma$ is
{\it regular with
valency $k$} whenever each vertex of $\Gamma$ is adjacent to exactly  $k$
distinct vertices of $\Gamma$. We say $\Gamma$ is {\it distance-regular}
whenever for all integers $h,i,j \ (0\le h,i,j \le D)$ and for all $x,y
\in X$ with $\partial(x,y)=h$, the number 
$$p^h_{ij}: = \left|\{z\in X|\ \partial(x,z) = i,\ \partial(y,z)=j\}\right| $$ is independent of $x$ and $y$. The
integers $p^h_{ij}$ are called the {\it intersection numbers} for
$\Gamma$. We abbreviate $a_i:=p_{1i}^i \hskip 5pt (0\le i\le D),\hskip 5pt
b_i:=p_{1i+1}^i \hskip 5pt (0\le i\le D-1),\hskip 5pt c_i:=p_{1i-1}^i
\hskip 5pt (1\le i\le D),$ and $ k_i:=p_{ii}^0 \hskip 5pt (0\le i\le D)$.
We observe $a_0=0$, $c_1=1$, and $k_0=1$. For notational convenience we
define $c_0:=0$ and $b_D:=0$. For the rest of this paper we assume $\Gamma$
is distance-regular with diameter $D\geq 3$. \\

We have a few comments.
The graph $\Gamma$ is regular with valency $k=k_1=b_0$. Moreover
\begin{equation}\label{eq:intnos}
c_i+a_i+b_i = k \qquad (0\le i \le D).
\end{equation}

\noindent It is known, by \cite[Chapter 3, Proposition 1.2]{bi}
\begin{equation}\label{eq:kifrac}
k_i = {{b_0b_1 \cdots b_{i-1}}\over{c_1c_2 \cdots c_i}} \qquad (0\le
i \le D).
\end{equation} 

\noindent
We remark $c_{i-1}\leq c_i$
and $b_{i-1}\geq b_i$
for $1 \leq i \leq D$ \cite[Proposition 4.1.6]{bcn}.\\

We recall the Bose-Mesner algebra of $\Gamma$. 
Let $\Mat_X({\C})$ denote the ${\C}$-algebra consisting of the matrices which have
rows and columns indexed by $X$ and entries in ${\C}$.  For $0 \le i \le D$ 
let $A_i$ denote the  matrix in $\Mat_X({\C})$ with entries
$${(A_i)_{xy} = \cases{1, & if $\;\partial(x,y) = i$\cr
0, & if $\;\partial(x,y) \ne i$\cr}} \qquad (x,y\in X).$$
We call $A_i$ the {\it ith distance matrix} of $\Gamma$.  We have
$A_0 = I,\
\sum_{i=0}^D A_i= J,\
A_i^t = A_i \ (0\le i \le D),\
A_iA_j = \sum_{h=0}^{D}p_{ij}^hA_{h}\ (0\le i,j \le D)$
where $J$ denotes the all 1's matrix.
The matrices $A_0, A_1, \ldots , A_D$ form a basis for a
commutative semi-simple ${\C}$-algebra $M$, called the {\it
Bose-Mesner
algebra of} $\Gamma$. By \cite[Section 2.3]{bi} $M$ has 
second basis $E_0, E_1, \ldots , E_D$ such that
$
E_0 = |X|^{-1}J,\
\sum_{i=0}^D E_i = I,\
E_i^t = E_i\ (0\le i \le D), \
E_iE_j = \delta_{ij}E_i \ (0\le i,j \le D).\
$
The  $E_0, E_1, \ldots , E_D$ are called the {\it primitive idempotents}
of $\Gamma$. We call $E_0$ the \emph{trivial} idempotent of $\Gamma$. \\

We set $A:=A_1$ and call this the  \emph{adjacency matrix} of $\Gamma$.
Let $\theta_0, \theta_1, \ldots,  \theta_D$ denote the complex scalars such that
$A = \sum_{i=0}^{D} \theta_iE_i.$
It is known  $\theta_0=k$, and that
$\theta_0, \theta_1, \ldots, \theta_D$ are
distinct real
numbers \cite[Chapter 3, Theorem 1.3]{bi}.
We refer to $\theta_i$ as the {\it eigenvalue of}
$\Gamma$ {\it associated with} $E_i \ (0 \le i \le D)$.
We call
$\theta_0$ the
{\it trivial eigenvalue of} $\Gamma$.
For each integer $i \ (0 \le i \le D)$, let
 $m_i$ denote the rank of $E_i$.  We refer to  $m_i$ as the {\it
multiplicity} of
$E_i\ ({\rm or\ }\theta_i).$
\\

Let $\theta$ denote an eigenvalue of $\Gamma$, let
$E$ denote the associated primitive idempotent, and let $m$
denote the multiplicity of $E$.  By \cite[Section 2.3]{bi}
there exists a sequence of real numbers $\sigma_0, \sigma_1,
\ldots, \sigma_D$ such
that
$$E = |X|^{-1}m\sum_{i=0}^D\sigma_iA_i.$$
It follows from \cite[Chapter 2, Proposition 3.3 (iii)]{bi} that
$\sigma_0=1$.
We call
$\sigma_0, \sigma_1, \ldots, \sigma_D$ the
{\it cosine sequence of} $\Gamma$ {\it associated with} $\theta$. 
We shall often abbreviate $\sigma_1$ by $\sigma$.  Let
$\sigma_0, \sigma_1, \ldots, \sigma_D, \theta$
denote real numbers.
By \cite[Proposition
4.1.1]{bcn} the following are equivalent: (i) $\theta$ is an eigenvalue of $\Gamma$ and
$\sigma_0, \sigma_1, \ldots, \sigma_D$ is the corresponding cosine sequence;
(ii) $\sigma_0 = 1$ and
$c_i\sigma_{i-1}+a_i\sigma_i+b_i\sigma_{i+1}=\theta\sigma_i
\ (0\le i \le D),$
where $\sigma_{-1}$ and $\sigma_{D+1}$ are indeterminates.\\

We end this section with a comment.
\begin{lemma}\label{lem:a1notzero}{\rm \cite[Proposition 5.5.1]{bcn}}
Let $\Gamma$ denote a distance-regular graph with diameter $D \ge
3$.   Suppose $a_1 \not=0$. Then $a_i \not=0$ for $1 \leq i \leq
D-1$.
\end{lemma}

\section{Pseudo cosine sequences}

We now recall a sequence of real numbers which generalizes the cosine sequence. 

\begin{definition}\label{def:pcs}\cite{jurter, terweng}
{\rm
Let $\Gamma$ denote a distance-regular graph with diameter $D \ge 3$.
For $\theta
\in \R$, by  the \emph{pseudo cosine sequence for
$\theta$} we mean the sequence of real numbers $\sigma_0,\sigma_1,
\dots, \sigma_D$
such that $\sigma_0=1$ and 
\begin{equation}\label{eq:recursion}
c_i\sigma_{i-1}+a_i\sigma_i+b_i\sigma_{i+1}=\theta\sigma_i \quad (0 \le i
\le D-1),
\end{equation} 
where $\sigma_{-1}$ is an indeterminate.
For notational
convenience we abbreviate $\sigma:=\sigma_1$.
We remark 
  $\theta$ is
determined by its pseudo cosine sequence; indeed
 $\theta=k \sigma$, where $k$ denotes the valency of $\Gamma$.
 }
\end{definition}
 
 We have some remarks on Definition \ref{def:pcs}.  Let $\Gamma$ denote
 a distance-regular graph with diameter $D\ge 3$. 
 Pick $\theta\in\R$ and let
$\sigma_0, \sigma_1, \dots, \sigma_D$
denote the corresponding pseudo cosine sequence.  Then $\theta$
is an eigenvalue of $\Gamma$ if and only if 
$c_D\sigma_{D-1}+a_D\sigma_D=\theta\sigma_D$.  In this case
$\sigma_0, \sigma_1, \dots, \sigma_D$ is the cosine sequence for $\theta$.
Now let $\sigma_0, \sigma_1, \dots, \sigma_D$ denote the pseudo cosine sequence
for $\theta=k$, where $k$ is the valency of $\Gamma$.  
Comparing
(\ref{eq:intnos}),
(\ref{eq:recursion}) we find
$\sigma_i=1 $ for $0 \leq i \leq D$.
By the {\it trivial} pseudo cosine sequence we mean
the pseudo cosine sequence for $k$. This sequence is a cosine
sequence since $k$ is an eigenvalue of $\Gamma$.

\begin{lemma}\label{lem:CD}{\bf (Christoffel-Darboux)}
Let $\Gamma$ denote a distance-regular graph with diameter
$D\ge 3$.  Let $\sigma_0, \sigma_1, \ldots, \sigma_D$
 and $\rho_0, \rho_1, \ldots, \rho_D$ denote pseudo cosine sequences of $\Gamma$.  Then
\begin{equation}\label{eq:CD}
({\sigma-\rho})\sum_{h=0}^{i}k_h\sigma_h\rho_h = {b_1b_2 \cdots b_{i}
\over{c_1c_2\cdots
c_{i}}}\bigl({\sigma_{i+1}\rho_{i}-\sigma_{i}\rho_{i+1}} \bigr)
\qquad\qquad (0\le i \le D-1). \end{equation} 
\end{lemma}

\begin{proof}
For $0 \leq h \leq i$ we have
\begin{eqnarray}
\sigma \sigma_h k &=& c_h\sigma_{h-1} + a_h \sigma_h + b_h \sigma_{h+1},
\label{eq:sigrec} \\ \rho \rho_h k&=& c_h\rho_{h-1} + a_h \rho_h + b_h
\rho_{h+1}. \label{eq:rhorec} 
\end{eqnarray} 
Evaluating the left-hand side
of (\ref{eq:CD}) using (\ref{eq:sigrec}), (\ref{eq:rhorec}), and
simplifying the result using (\ref{eq:kifrac}) we routinely obtain the
right-hand side of (\ref{eq:CD}). \end{proof}

\begin{lemma}\label{lem:recur}
Let $\Gamma$ denote a distance-regular graph
with diameter $D\ge 3$. Let
$\sigma_0, \sigma_1, \ldots, \sigma_D$ denote real numbers.
Then the following are equivalent.
\begin{enumerate}
\item $\sigma_0, \sigma_1, \ldots, \sigma_D$ is a pseudo cosine sequence. 
\item $\sigma_0 = 1$ and
\begin{equation}\label{eq:recur}
c_i(\sigma_{i-1}-\sigma_i)-
b_i(\sigma_i-\sigma_{i+1})=k(\sigma-1)\sigma_i \qquad\qquad (1\le i \le
D-1). \end{equation} 
\item $\sigma_0 = 1$ and
\begin{equation}\label{eq:CD2} ({\sigma-1})\sum_{h=0}^{i}k_h\sigma_h =
{b_1b_2 \cdots b_{i}\over{c_1c_2 \cdots
c_{i}}}\bigl({\sigma_{i+1}-\sigma_{i}}\bigr) \qquad\qquad (0\le i \le
D-1). \end{equation} \end{enumerate} \end{lemma} 
\begin{proof}
(i)$\Longrightarrow$(iii) To show
(\ref{eq:CD2}),  in Lemma \ref{lem:CD} let $\rho_0, \rho_1, \ldots, \rho_D$
denote the pseudo cosine sequence for $k$.\\
(iii)$\Longrightarrow$(ii) We show (\ref{eq:recur}) holds.
  Pick an integer
$i \ (1\le i \le D-1)$. By (\ref{eq:CD2}) (with $i$ replaced by $i-1)$ we
have  
\begin{equation}\label{eq:CD3}
(\sigma-1)\sum_{h=0}^{i-1}k_h\sigma_h =
{b_1b_2 \cdots b_{i-1}\over{c_1c_2 \cdots
c_{i-1}}}\bigl({\sigma_{i}-\sigma_{i-1}}\bigr). \end{equation}
 Subtracting equation (\ref{eq:CD3}) from
equation (\ref{eq:CD2}) and eliminating $k_i$ from the result using
(\ref{eq:kifrac}) we get (\ref{eq:recur}) as desired.\\
(ii)$\Longrightarrow$(i) We show
\begin{equation}\label{eq:recursion2}
c_i\sigma_{i-1}+a_i\sigma_i+b_i\sigma_{i+1}=k\sigma\sigma_i 
\end{equation}
for $0 \le i \le D-1$.
Clearly (\ref{eq:recursion2}) holds for  $i=0$. To show (\ref{eq:recursion2}) holds
for  $1 \le i \le D-1$
simplify (\ref{eq:recur}) using (\ref{eq:intnos}).  We now have (\ref{eq:recursion2}) for $0 \le i \le D-1$.  Applying Definition \ref{def:pcs} we find
$\sigma_0, \sigma_1, \ldots, \sigma_D$ is the pseudo cosine sequence for $k\sigma$.
 \end{proof}

We mention another characterization of the pseudo cosine
sequences.

\begin{lemma}\label{lem:recur2}
Let $\Gamma$ denote a distance-regular graph
with diameter $D\ge 3$. Let
$\sigma_0, \sigma_1, \ldots, \sigma_D$ denote real numbers.
Then the following are equivalent.
\begin{enumerate}
\item $\sigma_0, \sigma_1, \ldots, \sigma_D$ is a pseudo cosine sequence.
\item $\sigma_0 = 1$ and
\begin{equation}\label{eq:recurb2}
b_i(\sigma_{i-1}-\sigma_{i+1}) =
k(\sigma_{i-1}-\sigma \sigma_i)
-
a_i(\sigma_{i-1}-\sigma_i)
 \qquad\qquad (1\le i \le D-1).
\end{equation}
\item $\sigma_0 = 1$ and
\begin{equation}\label{eq:recurc2}
c_i(\sigma_{i+1}-\sigma_{i-1}) =
k(\sigma_{i+1}-\sigma \sigma_i)
-
a_i(\sigma_{i+1}-\sigma_i)
 \qquad\qquad (1\le i \le D-1).
\end{equation}
\end{enumerate}
\end{lemma}
\begin{proof}
(i)$\Longrightarrow$(ii) Lemma \ref{lem:recur}(i) holds. Applying that lemma
and eliminating $c_i$ in
(\ref{eq:recur})  using  (\ref{eq:intnos}) we get (\ref{eq:recurb2}).\\
(ii)$\Longrightarrow$(iii) Eliminating $b_i$ in (\ref{eq:recurb2}) using  (\ref{eq:intnos}) we get (\ref{eq:recurc2}).\\
(iii)$\Longrightarrow$(i) We show (\ref{eq:recur}) holds.  To do this we eliminate
$a_i$ in (\ref{eq:recurc2}) using (\ref{eq:intnos}).  It follows $\sigma_0, \sigma_1, \ldots, \sigma_D$ is a pseudo cosine sequence in view of 
Lemma \ref{lem:recur}(i),(ii).
\end{proof}

We finish this section with some results which we will find useful later in the paper.

\begin{lemma}
\label{lem:doublezero}
Let $\Gamma$ denote a distance-regular graph with diameter
$D\ge 3$. Let 
$\sigma_0, \sigma_1, \dots, \sigma_D$ denote a pseudo cosine
sequence of $\Gamma$. Then at least one of $\sigma_i, \sigma_{i+1}$
is nonzero for $0 \leq i \leq D-1$.
\end{lemma}
\begin{proof}
Assume there exists an integer $i$ $(0 \leq i \leq D-1)$ such
that $\sigma_i=0$ and $\sigma_{i+1}=0$. Recall $\sigma_0=1$ so
$i \geq 1$. Without loss of generality we may assume $\sigma_{i-1}\not=0$.
This is inconsistent with  (\ref{eq:recursion}) 
and the fact that $c_i\ne 0$.  The result follows. 
\end{proof}

\begin{lemma}
\label{lem:triplemult}
Let $\Gamma$ denote a distance-regular graph with diameter
$D\ge 3$. Let
$\sigma_0, \sigma_1, \dots, \sigma_D$ denote a nontrivial
pseudo cosine
sequence of $\Gamma$. Then for $1 \leq i \leq D-1,$
at least one of
 $\sigma_{i-1}-\sigma_{i}$,
 $\sigma_{i}-\sigma_{i+1}$
is nonzero. 
\end{lemma}
\begin{proof}
Suppose
there exists an integer $i$
$(1 \leq i \leq D-1)$
such that each of
 $\sigma_{i-1}-\sigma_{i}$,
 $\sigma_{i}-\sigma_{i+1}$   is zero.
Then
 $\sigma_{i-1}, \sigma_{i},
 \sigma_{i+1}$
coincide and
this common value is nonzero by
Lemma
\ref{lem:doublezero}.
Evaluating
(\ref{eq:recur})
using this
we find $\sigma=1$. Now
$\sigma_0, \sigma_1, \dots, \sigma_D$ is the trivial
pseudo cosine sequence, a contradiction.
\end{proof}

\begin{lemma}
\label{lem:grandmult}
Let $\Gamma$ denote a 
distance-regular graph with diameter
$D\ge 3$ and valency $k$.  Let
$\sigma_0, \sigma_1, \dots, \sigma_D$ denote the
pseudo cosine
sequence of $\Gamma$ associated with $-k$.
 Suppose there exists an integer $i$ $(1 \leq i \leq D-2)$
 such that $\sigma_{i-1}+\sigma_{i}=0$ and
 $\sigma_{i+1}+\sigma_{i+2}=0$.
Then  $a_i=0$, $a_{i+1}=0$, and
 $\sigma_{i}+\sigma_{i+1}=0$.
\end{lemma}
\begin{proof}
Observe $\sigma_i\not=0$; otherwise each of
$\sigma_{i-1},
\sigma_i$ is zero, contradicting
Lemma
\ref{lem:doublezero}.
Also $\sigma_{i+1}\not=0$; otherwise each of
$\sigma_{i+1},
\sigma_{i+2}$ is zero, contradicting
Lemma
\ref{lem:doublezero}.
Applying
(\ref{eq:recursion}) at $i$ and $i+1$ with $\theta=-k$
we find
both
$b_i\sigma_{i+1}= \sigma_i(c_i-a_i-k)$
and
$c_{i+1}\sigma_{i}= \sigma_{i+1}(b_{i+1}-a_{i+1}-k)$.
Combining these two equations and using
(\ref{eq:intnos}) we routinely obtain
\begin{eqnarray}
(k+a_i-c_i)a_{i+1} = -a_ic_{i+1}.
\label{almostthere}
\end{eqnarray}
In
(\ref{almostthere}) the left-hand side is nonnegative
and the right-hand side is nonpositive. Therefore both
sides are zero. It follows each of $a_i, a_{i+1}$ is zero.
Setting $a_i=0$ in the equation
$b_i\sigma_{i+1}= \sigma_i(c_i-a_i-k)$
we routinely obtain $\sigma_i+\sigma_{i+1}=0$.
\end{proof}

\section{The tight relation}

\begin{definition}
\label{def:tight}
{\rm
Let $\Gamma$ denote a distance-regular graph with diameter
$D\ge 3$.
We introduce a binary symmetric relation on the set of all
pseudo cosine sequences of $\Gamma$. We call this the \emph{tight}
relation. Let
$\sigma_0, \sigma_1, \ldots, \sigma_D$  and
$\rho_0, \rho_1, \ldots, \rho_D$ denote pseudo cosine sequences
of $\Gamma$. We  say this pair of
sequences is \emph{tight} whenever
 $\sigma_0\rho_0, \sigma_1\rho_1, \ldots, \sigma_D\rho_D$ is a
pseudo cosine sequence of $\Gamma$.}
 \end{definition}

\begin{definition}
\label{def:tightp}
{\rm Let $\Gamma$ denote a distance-regular graph with diameter
$D\ge 3$.
We define a binary symmetric relation on
$\R$
which we call the \emph{tight}
relation.
For
$\theta, \theta' \in \R$,
we say the pair
$\theta, \theta'$ is \emph{tight} whenever the
corresponding pseudo cosine sequences form a tight pair.}
\end{definition}

\noindent
Let $\Gamma$ denote a distance-regular graph with diameter
$D\ge 3$ and valency $k$.
We give an example of a  tight pair of
pseudo cosine sequences.
Let $\sigma_0, \sigma_1, \dots, \sigma_D$
denote the pseudo cosine sequence for $k$
and recall $\sigma_i=1 $ for $0 \leq i \leq D$.
Let $\rho_0, \rho_1, \ldots, \rho_D$ denote any pseudo cosine sequence.
Of course $\sigma_i \rho_i=\rho_i$ for $0 \leq i \leq D$. Applying
Definition
\ref{def:tight} we find
$\sigma_0, \sigma_1, \dots, \sigma_D$
and
$\rho_0, \rho_1, \ldots, \rho_D$ form a tight pair.
Consequently we have the following.

\begin{example}
\label{ex:triv}
{\rm
Let $\Gamma$ denote a distance-regular graph with diameter
$D\geq 3$ and valency $k$.
Then for
all $\theta \in \R$ the pair
$\theta, k$ is tight.}
\end{example}

\begin{lemma}
\label{lem:tpev}
Let $\Gamma$ denote a distance-regular graph with diameter
$D\ge 3$ and valency $k$.
Let $\theta, \theta'$ denote real numbers other than $k$,
and assume
$\theta, \theta'$ form a tight pair.
Then
\begin{equation}
\Biggl(\theta + \frac{k}{a_1+1}\Biggr)
\Biggl(\theta' + \frac{k}{a_1+1}\Biggr) = -\frac{ka_1b_1}{(a_1+1)^2}.
\label{eq:ththp}
\end{equation}
\end{lemma}
\begin{proof} Let
 $\sigma_0, \sigma_1, \ldots, \sigma_D$
(resp.
$\rho_0, \rho_1, \ldots, \rho_D$) denote the pseudo cosine
sequence for $\theta$ (resp. $\theta'$).
Define $\gamma_i:=\sigma_i\rho_i$ for $0 \leq i \leq D$. Observe
$\gamma_0, \gamma_1, \ldots, \gamma_D$
is the pseudo cosine sequence for $\psi$, where
$\psi:=k \sigma \rho$.
Setting $i=1$ in
(\ref{eq:recursion}) and using $\theta=k\sigma$ we find
\begin{equation}
1 + a_1\sigma + b_1 \sigma_2 = k \sigma^2.
\label{eq:oneof3}
\end{equation}
Similarly
\begin{eqnarray}
1 + a_1\rho + b_1 \rho_2 &=& k \rho^2,
\label{eq:twoof3}
\\
1 + a_1\sigma \rho + b_1 \sigma_2\rho_2 &=& k \sigma^2\rho^2.
\label{eq:threeof3}
\end{eqnarray}
To obtain
(\ref{eq:ththp}),
first solve
 (\ref{eq:oneof3}) and (\ref{eq:twoof3}) for
$\sigma_2 $ and $\rho_2$ respectively.
Then evaluate
(\ref{eq:threeof3}) using this
and simplify the result using
 $\theta =k \sigma$, $\theta'=k\rho$ and (\ref{eq:intnos}).
\end{proof}

\noindent In this paper, we will find all the tight pairs.
We will break down the argument into the following cases 
(i) $a_i= 0$ for $0 \leq i \leq D-1$; (ii)
$a_i= 0$ for $0 \leq i \leq D-2$ and $a_{D-1}\not=0$;
(iii) $a_1=0$ and there exists $i$ $(2 \leq i \leq D-2)$ such
that $a_i\not=0$;
(iv) $a_1 \not=0$.

\section{The case $a_i=0 \;(0 \leq i\leq D-1)$
}

Let $\Gamma$ denote a distance-regular graph with diameter
$D\ge 3$. Recall $\Gamma$ is {\it bipartite}
whenever $a_i=0$ for $0 \leq i \leq D$.
We say $\Gamma$ is {\it almost bipartite}
whenever $a_D\not=0$ and $a_i=0$ for $0 \leq i \leq D-1$.

\begin{lemma}
\label{lem:bipdual}
Let $\Gamma$ denote a
distance-regular graph with diameter
$D\ge 3$.  Assume $\Gamma$ is bipartite or almost bipartite.
Given $\theta \in \R$, let
$\sigma_0, \sigma_1, \ldots, \sigma_D$ denote
the pseudo cosine sequence for $\theta$ and let
$\rho_0, \rho_1, \ldots, \rho_D$ denote
 the pseudo cosine sequence for $-\theta$.
Then $\rho_i=(-1)^i\sigma_i$ for $0 \leq i \leq D$.
\end{lemma}

\begin{proof}
Observe the desired result holds for $i=0$ and $i=1$ since  $\sigma_0=1=\rho_0$
and $\sigma=\theta k^{-1}=-\rho$. To show the result holds for $2\le i \le D$ we use induction and the fact that $c_{i-1}\sigma_{i-2}+b_{i-1}\sigma_{i}=\theta\sigma_{i-1}$ and 
$c_{i-1}\rho_{i-2}+b_{i-1}\rho_{i}=-\theta\rho_{i-1}$.
\end{proof}

\begin{corollary}
\label{cor:triv2}
Let $\Gamma$ denote a
distance-regular graph with diameter
$D\ge 3$ and valency $k$.
 Assume $\Gamma$ is bipartite or almost bipartite.
Let $\rho_0, \rho_1, \ldots, \rho_D$ denote
the pseudo cosine sequence for $-k$. Then
$\rho_i=(-1)^i$ for $0 \leq i \leq D$.
\end{corollary}
\begin{proof} Referring to Lemma
 \ref{lem:bipdual},
set $\theta=k$
and recall $\sigma_i=1 $ for $0 \leq i \leq D$.
\end{proof}

\begin{theorem}
\label{thm:bpclass}
Let $\Gamma$ denote a
distance-regular graph with diameter
$D\ge 3$ and valency $k$.
 Assume $\Gamma$ is bipartite or almost bipartite.
\begin{enumerate}
\item The pair $\theta, k$ is tight for all $\theta \in \R$.
\item The pair $\theta, -k$ is tight for all $\theta \in \R$.
\item $\Gamma$ has no further tight pairs.
\end{enumerate}
\end{theorem}
\begin{proof}(i)  This is from Example
\ref{ex:triv}.
\\
\noindent (ii) Let $\sigma_0, \sigma_1, \ldots, \sigma_D$ denote
the pseudo cosine sequence for $\theta$, and let  
$\rho_0, \rho_1, \ldots, \rho_D$ denote
the pseudo cosine sequence for $-k$.
We show the sequences $\sigma_0, \sigma_1, \ldots, \sigma_D$
and $\rho_0, \rho_1, \ldots, \rho_D$  form a tight pair.
For $0 \le i \le D$ define
$\gamma_i:= \sigma_i\rho_i$.  Recall $\rho_i=(-1)^i$ for $0 \le i \le D$ by Corollary \ref{cor:triv2} so $\gamma_i=(-1)^i\sigma_i$.
 Observe $\gamma_0, \gamma_1, \ldots, \gamma_D$ 
 is the pseudo cosine sequence  for $-\theta$ by Lemma  \ref{lem:bipdual}.
 We have now shown the sequences $\sigma_0, \sigma_1, \ldots, \sigma_D$
and $\rho_0, \rho_1, \ldots, \rho_D$  form a tight pair.
 It follows the pair $\theta, -k$ is tight.
\\
\noindent (iii)
For $\theta\in\R,\ \theta'\in\R$,
 suppose
 $\theta, \theta'$ is a tight pair but
$\theta \not=k$, $\theta' \not=k$.
We show $\theta=-k$ or
$\theta'=-k$.
Applying Lemma
\ref{lem:tpev}
we find $\theta, \theta'$ satisfy
(\ref{eq:ththp}). By this and since  $a_1=0$
we find $(\theta+k)(\theta'+k)=0$.
Now $\theta=-k$ or
$\theta'=-k$.
The result follows.
\end{proof}

\section{Tight pairs $\theta, \theta'$ with $\theta=\theta'$.}

Let $\Gamma$ denote a
distance-regular graph with diameter
$D\ge 3$ and valency $k$.
Let $\theta, \theta'$ denote real numbers other than $k$,
and assume the pair
$\theta, \theta'$ is tight.
Conceivably
$\theta = \theta'$;
we consider when this occurs.
First suppose $\Gamma$ is bipartite or almost bipartite.
>From Lemma
\ref{thm:bpclass} we find
$\theta=\theta'$ if and only if
each of
$\theta, \theta'$ is equal to $-k$.
We now treat the case in which $\Gamma$ is neither bipartite nor
almost bipartite.

\begin{lemma}
\label{lem:dist}
Let $\Gamma$ denote a distance-regular graph with diameter
$D\ge 3$ and valency $k$. Assume $\Gamma$ is neither bipartite nor almost bipartite.
Let $\theta, \theta'$ denote real numbers other than $k$,
and assume the pair
$\theta, \theta'$ is tight.
Then
 $\theta'\not=\theta$.
\end{lemma}
\begin{proof}
Applying Lemma
\ref{lem:tpev}
we find  $\theta, \theta'$ satisfies
(\ref{eq:ththp}).
Suppose $\theta=\theta'$.
Then
the left-hand side of (\ref{eq:ththp}) is
a square so it is nonnegative.
The right-hand side of (\ref{eq:ththp}) is nonpositive so both
sides  of (\ref{eq:ththp}) are equal to $0$. Now $a_1=0$. Moreover
$\theta$ and $\theta'$ are equal to $-k$.
Let $\sigma_0, \sigma_1, \ldots, \sigma_D$ denote the pseudo cosine
sequence for $\theta$.
Since the pair $\theta ,\theta$ is tight
we find
$\sigma^2_0, \sigma^2_1, \ldots, \sigma^2_D$ is a pseudo cosine
sequence. This is a pseudo cosine sequence for $\psi$ where
$\psi: = k \sigma^2$.
Recall $\sigma = \theta/k$ and we mentioned
$\theta=-k$ so
$\sigma=-1$. Apparently $\psi=k$
so $\sigma^2_i=1$ for $0 \leq i \leq D$.
We assume $\Gamma$ is neither bipartite nor almost bipartite.
Therefore there exists
an integer $r$ $(0 \leq r\leq D-1)$ such that $a_r\not=0$.
Assume we have picked the minimal such $r$.
Observe $r\geq 2$ since $a_0=0$ and $a_1=0$.
For $0 \leq i < r$ we
apply
(\ref{eq:recursion})
with $\theta=-k$
and
$a_i=0$ to obtain
$\sigma_j=(-1)^j$ for $0 \leq j \leq r$.
We now
apply
(\ref{eq:recursion})
with $i=r$ and
$\theta=-k$.
 Simplifying the result using
$\sigma_{r-1}=(-1)^{r-1}$,
$\sigma_r=(-1)^r$
and (\ref{eq:recurb2})
we obtain
\begin{equation}
1-(-1)^{r+1}\sigma_{r+1}=-2a_r/b_r.
\label{eq:getcont}
\end{equation}
>From our above comments
$\sigma_{r+1}^2=1$ so
 $\sigma_{r+1}=1$ or
 $\sigma_{r+1}=-1$. It follows 
the left-hand side of
(\ref{eq:getcont}) is nonnegative. By construction
$a_r>0$ so the right-hand side
of
(\ref{eq:getcont}) is negative. We now have a contradiction.
We conclude $\theta\not=\theta'$.
\end{proof}

\section{The auxiliary parameter}

We now introduce a parameter which will help us describe the tight pairs.

\begin{theorem}
\label{thm:main1}
Let $\Gamma$ denote a distance-regular graph with diameter
$D\ge 3$. Let
$\sigma_0, \sigma_1, \ldots, \sigma_D$  and
$\rho_0, \rho_1, \ldots, \rho_D$ denote
nontrivial pseudo cosine sequences
of $\Gamma$. Then the following are equivalent.
\begin{enumerate}
\item
$\sigma_0, \sigma_1, \ldots, \sigma_D$  and
$\rho_0, \rho_1, \ldots, \rho_D$
 form a tight pair.
\item There exists a real number $\varepsilon$ such that
\begin{eqnarray}
\sigma_i \rho_i - \sigma_{i-1} \rho_{i-1} =
\varepsilon(
\sigma_{i-1} \rho_i - \sigma_{i} \rho_{i-1}) \qquad \qquad
(1 \leq i \leq D).
\label{eq:epseq}
\end{eqnarray}
\end{enumerate}
\end{theorem}
\begin{proof}  Set
\begin{equation}\label{eq:gammai}
\gamma_i=\sigma_i\rho_i \ (0 \le i \le D).
\end{equation}
(i)$\Longrightarrow$(ii)  By Definition \ref{def:tight}
the sequence $\gamma_0, \gamma_1, \dots, \gamma_D$ is a pseudo cosine
sequence. First assume $\sigma \ne \rho$, and let
\begin{equation}\label{eq:epsilon} \varepsilon
={{\sigma\rho-1}\over{\rho-\sigma}}.
 \end{equation}
Pick an integer $i \ (1 \le i \le D)$.
By  Lemma \ref{lem:CD}, Lemma \ref{lem:recur} and 
(\ref{eq:epsilon}) we have
\begin{eqnarray}
\sigma_i\rho_i - \sigma_{i-1}\rho_{i-1} &=& \gamma_i - \gamma_{i-1}
      \nonumber \\
{ } &=& (\gamma -1)\
{{c_1c_2 \cdots c_{i-1}}\over{b_1b_2 \cdots b_{i-1}}}\
      \sum_{h=0}^{i-1}k_h\gamma_h
       \nonumber \\
{ } &=& (\sigma\rho - 1)\
{{c_1c_2 \cdots c_{i-1}}\over{b_1b_2 \cdots b_{i-1}}}\
      \sum_{h=0}^{i-1}k_h\sigma_h\rho_h
       \nonumber \\
{ } &=& (\sigma\rho - 1) \,
\frac{\sigma_i\rho_{i-1}-\sigma_{i-1}\rho_i}{\sigma-\rho}
\nonumber \\
{ } &=& \varepsilon(\sigma_{i-1}\rho_i - \sigma_i\rho_{i-1}),\nonumber
\end{eqnarray} as desired.
Next assume $\sigma= \rho$. Observe
$\Gamma$ is bipartite or almost bipartite in view of Lemma \ref{lem:dist}. 
By Theorem
\ref{thm:bpclass} and since $k\sigma, k\rho$ form a tight pair we
find each of $k\sigma, k \rho$ is equal to $-k$.  Now by Corollary
\ref{cor:triv2} each of $\sigma_i, \rho_i$ is equal to $(-1)^i$
for $0 \le i \le D$.
It follows (\ref{eq:epseq}) holds for any real number $\varepsilon$. \\
\noindent (ii)$\Longrightarrow$(i) Assume (\ref{eq:epseq}) holds for some
real number $\varepsilon$.  We show  $\gamma_0, \gamma_1, \dots,
\gamma_D$ is a pseudo cosine sequence. By Lemma \ref{lem:recur}(i),(iii)
it suffices to show
\begin{equation}\label{eq:CD3p}
({\gamma-1})\sum_{h=0}^{i}k_h\gamma_h = {b_1b_2 \cdots b_{i}\over{c_1c_2
\cdots c_{i}}}\bigl({\gamma_{i+1}-\gamma_{i}}\bigr)
 \end{equation}
for $0\le i \le D-1.$
Pick an integer $i \ (0 \le i \le D-1)$.
By Lemma \ref{lem:CD} we have
\begin{equation}\label{eq:CD4}
({\sigma-\rho})\sum_{h=0}^{i}k_h\sigma_h\rho_h = {b_1b_2 \cdots b_{i}
\over{c_1c_2\cdots
c_{i}}}\bigl({\sigma_{i+1}\rho_{i}-\sigma_{i}\rho_{i+1}} \bigr).
\end{equation} 
 Setting $i=1$ in
(\ref{eq:epseq}) we get
\begin{equation}\label{eq:epseqat1}
\sigma\rho-1=\varepsilon(\rho-\sigma). \end{equation} Multiplying equation
(\ref{eq:CD4}) by $\varepsilon$ and simplifying the result using
(\ref{eq:epseqat1}), (\ref{eq:epseq}) and (\ref{eq:gammai}) gives
(\ref{eq:CD3p}) as desired.

\end{proof}

\begin{definition}\label{def:auxpar}
{\rm Let $\Gamma$ denote a distance-regular graph with diameter $D\ge
3$. Suppose we are given nontrivial pseudo cosine sequences
 $\sigma_0, \sigma_1, \ldots, \sigma_D$
and $\rho_0, \rho_1, \ldots, \rho_D$ which form a tight pair. By an
\emph{auxiliary parameter} for this pair
we mean a real number  $\varepsilon$ which satisfies  
Theorem \ref{thm:main1}(ii).
We comment on the uniqueness of the auxiliary parameter.
Suppose $\sigma=\rho$. Then $\varepsilon$ is an auxiliary
parameter for all $\varepsilon \in \R$. Suppose
$\sigma\not=\rho$. Setting $i=1$ in (\ref{eq:epseq})
we find the auxiliary parameter satisfies
\begin{eqnarray}\label{eq:auxform}
\varepsilon = \frac{\sigma \rho -1}{\rho - \sigma}.
\end{eqnarray}
In this case the auxiliary parameter is unique.}
\end{definition}

\section{When can $\sigma_{i-1}= \sigma_i$?}

\indent Let\  $\Gamma$\ denote a distance-regular graph with
diameter $D\geq 3$.\ \ Assume $\sigma_0, \sigma_1, \ldots, \sigma_D$ and
 $\rho_0, \rho_1, \ldots, \rho_D$ denote nontrivial
pseudo cosine sequences which form a tight pair. Let $\varepsilon$ denote a corresponding
auxiliary parameter.  A bit later  in this paper,
we discuss certain rational expressions involving $\sigma_0,
\sigma_1, \ldots, \sigma_D$  and $\varepsilon$
which contain in the denominator the
factors $\sigma_{i-1}-\sigma_i$  for $1
\leq i \leq D$. To prepare for this discussion, we investigate how
these factors can be zero.

\begin{lemma}
\label{lem:epsone} Let $\Gamma$ denote a distance-regular graph
with diameter $D\ge 3$.  Suppose we are given nontrivial
pseudo cosine sequences $\sigma_0, \sigma_1, \ldots, \sigma_D$ and
$\rho_0, \rho_1, \ldots, \rho_D$ which form a tight pair. Let us
assume this pair has auxiliary parameter $\varepsilon=1$. Then
\begin{equation}
(\sigma_{i-1}-\sigma_i)(\rho_{i-1}+\rho_i)=0 \qquad \qquad  (1 \leq i \leq D).
\label{eq:oneiszero}
\end{equation}
\end{lemma}
\begin{proof} To get
(\ref{eq:oneiszero}) set $\varepsilon=1$ in
(\ref{eq:epseq}).
\end{proof}

\begin{lemma}
\label{lem:threeres} Let $\Gamma$ denote a distance-regular graph
with diameter $D\ge 3$. Let us assume $\Gamma$ is neither
bipartite nor almost bipartite. Suppose we are given nontrivial
pseudo cosine sequences $\sigma_0, \sigma_1, \ldots, \sigma_D$ and
$\rho_0, \rho_1, \ldots, \rho_D$ which form a tight pair. Let us
assume this pair has auxiliary parameter $\varepsilon=1$. Then
the following (i)--(iii) hold.
\begin{enumerate}
\item $\rho_i = (-1)^i$ for $0 \leq i \leq D-1$ and $\rho_D \not=(-1)^D$.
\item $\sigma_{D-1}=\sigma_D$.
\item $a_i=0$ for $0 \leq i \leq D-2$ and $a_{D-1}\not=0$.
\end{enumerate}
\end{lemma}
\begin{proof}
(i)
Recall $\rho_0=1$.
 Setting $i=1$ in
(\ref{eq:oneiszero}) and recalling $\sigma_0=1$, $\sigma\not=1$ we
find $\rho=-1$. Let $r$ denote  the maximal integer $(0 \leq r
\leq D)$ such that $\rho_i=(-1)^i$ for $0\leq i \leq r$. We show
$r=D-1$. From our preliminary comments $r \geq 1$. Observe $r
\not=D$; otherwise $\Gamma$ is bipartite or almost bipartite in
view of  Lemma \ref{lem:grandmult}. Suppose $r \leq D-2$. By the
construction $\rho_{r-1} + \rho_{r} =0$ and $\rho_r + \rho_{r+1}
\not=0$. Applying Lemma \ref{lem:epsone} we find
$\sigma_r=\sigma_{r+1} $. Observe
 $\sigma_{r+1}\not=\sigma_{r+2} $ by
Lemma
\ref{lem:triplemult} so
 $\rho_{r+1}+\rho_{r+2}=0 $ by
Lemma
\ref{lem:epsone}.
Now
 $\rho_{r}+\rho_{r+1}=0 $ by
Lemma
\ref{lem:grandmult}, a contradiction.
Hence $r=D-1$.
\\
\noindent (ii)
Apply
(\ref{eq:oneiszero})
with $i=D$  and observe  $\rho_{D-1} +\rho_D\not=0$ by (i) above.
\\
\noindent (iii) Recall $a_0=0$.  We find $a_i=0$ for $1 \le i \le
D-2$ by Lemma \ref{lem:grandmult} and since $\rho_i=(-1)^i$ for $0
\le i \le D-1$.  Observe $a_{D-1}\not=0$; otherwise $\Gamma$ is
bipartite or almost bipartite, contradicting our assumption.
\end{proof}

\begin{lemma}
Let $\Gamma$ denote a distance-regular graph with diameter $D\ge
3$.  Suppose we are given nontrivial pseudo cosine sequences
$\sigma_0, \sigma_1, \ldots, \sigma_D$ and $\rho_0, \rho_1, \ldots,
\rho_D$ which form a tight pair. Let us assume this pair has
auxiliary parameter $\varepsilon=-1$. Then
\begin{equation}
(\sigma_{i-1}+\sigma_i)(\rho_{i-1}-\rho_i)=0 \qquad \qquad  (1 \leq i \leq D).
\label{eq:oneiszero2}
\end{equation}
\end{lemma}
\begin{proof} To get
(\ref{eq:oneiszero2}) set $\varepsilon=-1$ in
(\ref{eq:epseq}).
\end{proof}

\begin{lemma}
\label{lem:threeres2} Let $\Gamma$ denote a distance-regular graph
with diameter $D\ge 3$. Let us assume $\Gamma$ is neither
bipartite nor almost bipartite. Suppose we are given nontrivial
pseudo cosine sequences $\sigma_0, \sigma_1, \ldots, \sigma_D$ and
$\rho_0, \rho_1, \ldots, \rho_D$ which form a tight pair. Let us
assume this pair has auxiliary parameter $\varepsilon=-1$. Then
the following (i)--(iii) hold.
\begin{enumerate}
\item $\sigma_i = (-1)^i$ for $0 \leq i \leq D-1$ and $\sigma_D \not=(-1)^D$.
\item $\rho_{D-1}=\rho_D$.
\item $a_i=0$ for $0 \leq i \leq D-2$ and $a_{D-1}\not=0$.
\end{enumerate}
\end{lemma}
\begin{proof}
In the proof of
Lemma
\ref{lem:threeres} interchange the roles of
$\sigma_0, \sigma_1, \ldots, \sigma_D$
and $\rho_0, \rho_1, \ldots, \rho_D$.
\end{proof}

\begin{lemma}
\label{lem:nodup} Let $\Gamma$ denote a distance-regular graph
with diameter $D\ge 3$. Suppose we are given nontrivial pseudo
cosine sequences $\sigma_0, \sigma_1, \ldots, \sigma_D$ and $\rho_0,
\rho_1, \ldots, \rho_D$ which form a tight pair. Let $\varepsilon $
denote a corresponding auxiliary parameter and assume
$\varepsilon \not\in \lbrace 1,-1\rbrace $. Then
$\sigma_{i-1}\not=\sigma_i$ and $\rho_{i-1}\not=\rho_i$ for $1
\leq i \leq D$.
\end{lemma}
\begin{proof}
First assume $\sigma=\rho$.  Observe 
$\Gamma$ is bipartite or almost bipartite in view of Lemma \ref{lem:dist}. 
Applying Theorem \ref{thm:bpclass} we find
  $\sigma_0, \sigma_1, \ldots, \sigma_D$
and $\rho_0, \rho_1, \ldots, \rho_D$ are pseudo
cosine sequences for $-k$. By Corollary \ref{cor:triv2}
we have $\sigma_i=\rho_i=(-1)^i$ for $ 0 \le i \le D$; so
$\sigma_{i-1}\not=\sigma_i$ and $\rho_{i-1}\not=\rho_i$ for $1
\leq i \leq D$ as desired. \\
Next assume $\sigma\not=\rho$.  Then
$\varepsilon$ is as in (\ref{eq:auxform}).
Suppose there exists an integer $i$ $(1
\leq i \leq D)$ such that $\sigma_{i-1}=\sigma_i$. Since $\sigma_0=1$ and
$\sigma\not=1$ we find $i \geq 2$. Without loss of generality we may
assume $\sigma_{i-2}\not=\sigma_{i-1}$. By Lemma \ref{lem:doublezero} the
common value of $\sigma_{i-1}, \sigma_i$ is nonzero. Setting
$\sigma_{i-1}=\sigma_i$ in (\ref{eq:epseq}) and simplifying we find
$(\varepsilon-1)(\rho_{i-1}-\rho_i)=0.$ We assume $\varepsilon \not=1$ so
$\rho_{i-1}=\rho_i$. From Lemma \ref{lem:recur}(ii) (with $i$ replaced by
$i-1$)  we find \begin{eqnarray}
c_{i-1}(\sigma_{i-2}-\sigma_{i-1})&=&k(\sigma -1)\sigma_{i-1},
\label{eq:cc} \\ c_{i-1}(\rho_{i-2}-\rho_{i-1})&=&k(\rho -1)\rho_{i-1}.
\label{eq:cc1} \end{eqnarray} Combining (\ref{eq:cc}), (\ref{eq:cc1}) we
obtain \begin{eqnarray}
(\sigma_{i-2}-\sigma_{i-1}) (\rho-1)\rho_{i-1} = (\rho_{i-2}-\rho_{i-1})
(\sigma-1)\sigma_{i-1}. \label{eq:frac} \end{eqnarray} By (\ref{eq:epseq})
(with $i$ replaced by $i-1$)  we find \begin{eqnarray}
\sigma_{i-1}\rho_{i-1}-\sigma_{i-2}\rho_{i-2} = \varepsilon (
\sigma_{i-2}\rho_{i-1}-\sigma_{i-1}\rho_{i-2}). \label{eq:iadj}
\end{eqnarray} Adding $(\sigma-1)\sigma_{i-1}$ times (\ref{eq:iadj}) and
$\varepsilon \sigma_{i-1}-\sigma_{i-2}$ times (\ref{eq:frac}), and
simplifying the result using
(\ref{eq:auxform}),
 we routinely find $\sigma_{i-2}=\sigma
\sigma_{i-1}$. Evaluating (\ref{eq:cc}) using this we find $c_{i-1}=k$, a
contradiction. \end{proof}

\begin{corollary}
\label{cor:epsgood}
Let $\Gamma$ denote a distance-regular graph with diameter
$D\ge 3$ and intersection number $a_1 \not=0$.
Suppose we are given nontrivial
pseudo cosine sequences $\sigma_0, \sigma_1, \ldots, \sigma_D$ and
$\rho_0, \rho_1, \ldots, \rho_D$ which form a tight pair.
Let $\varepsilon $ denote the corresponding
auxiliary parameter.
Then
 $\varepsilon  \not\in \lbrace 1, -1 \rbrace$.
Moreover
$\sigma_{i-1}\not=\sigma_i$ and $\rho_{i-1}\not=\rho_i$ for $1 \leq i
\leq D$.
\end{corollary}
\begin{proof}
We assume $a_1 \not=0$; therefore $\Gamma$ is neither bipartite
nor almost bipartite. Suppose $\varepsilon = 1$. By Lemma
\ref{lem:threeres}(iii) and since $D\geq 3$ we find $a_1 =0$ for a
contradiction.
 Suppose $\varepsilon = -1$.
By Lemma
\ref{lem:threeres2}(iii) and  since $D\geq 3$ we find
$a_1 =0$ for a contradiction.
Apparently $\varepsilon  \not\in \lbrace 1, -1 \rbrace$.
Applying
Lemma \ref{lem:nodup}
we find
$\sigma_{i-1}\not=\sigma_i$ and $\rho_{i-1}\not=\rho_i$ for $1 \leq i
\leq D$.
\end{proof}

\section{The case $a_i = 0 \;(0 \leq i \leq D-2)$, $a_{D-1}\not=0$.}

\begin{lemma}
\label{lem:pcformk}
Let
$\Gamma$ denote a distance-regular graph with diameter
$D\ge 3$ and valency $k$.
Assume
 $a_i = 0$ for $0 \leq i \leq D-2$ and $a_{D-1}\not=0$.
Let
 $\sigma_0, \sigma_1, \ldots, \sigma_D$
denote the  pseudo cosine sequence for $-k$.
Then $\sigma_i=(-1)^i $ for $0 \leq i\leq D-1$
and
\begin{eqnarray}
(-1)^D\sigma_D =  1 + 2a_{D-1}/b_{D-1}.
\label{eq:rhoD}
\end{eqnarray}
\end{lemma}
\begin{proof} 
The sequence  $\sigma_0, \sigma_1, \ldots, \sigma_D$
satisfies (\ref{eq:recursion}) with  $\theta=-k$.  It follows
$c_i\sigma_{i-1}+b_i\sigma_{i+1}=-k\sigma_i$ for $0 \le i \le D-2.$
Using this and (\ref{eq:intnos}), we routinely verify by induction that
$\sigma_i=(-1)^i$ for $0 \leq i\leq D-1.$
To get  (\ref{eq:rhoD}) evaluate
(\ref{eq:recursion}) at $i=D-1$ and simplify using (\ref{eq:intnos}).
\end{proof}

\noindent
Let $\Gamma$ denote a distance-regular graph with diameter
$D\ge 3$.
Assume
 $a_i = 0$ for $0 \leq i \leq D-2$ and $a_{D-1}\not=0$.
In what follows we will be concerned with the matrix
\begin{eqnarray}
\left(\begin{array}{cccccc}
0 & b_0 & & & & {\bf 0}\\
c_1 & 0 & b_1 & & & \\
 & c_2 & 0 & \cdot & & \\
 & & \cdot & \cdot &  b_{D-3}& \\
 & & & c_{D-2} & 0 & b_{D-2} \\
 {\bf 0} & & & & c_{D-1} & k-c_{D-1}
\end{array} \right).
\label{eq:mat}
\end{eqnarray}

\noindent We make a routine observation.

\begin{lemma}
\label{lem:Gmat}
Let $\Gamma$ denote a distance-regular graph with diameter
$D\ge 3$. Assume $a_i = 0$ for $0 \leq i \leq D-2$ and $a_{D-1}\not=0$.
Given scalars $\sigma_0,\sigma_1, \ldots, \sigma_{D-1}, \theta \in \R$,
the following are equivalent.
\begin{enumerate}
\item $\theta$ is an eigenvalue of the matrix (\ref{eq:mat}), and
 $(\sigma_0,\sigma_1, \ldots, \sigma_{D-1})^t$ is a corresponding
 eigenvector, scaled so that $\sigma_0=1$.
\item
 $\sigma_0,\sigma_1, \ldots, \sigma_{D-1}, \sigma_{D-1}$
is the pseudo cosine sequence for
 $\theta$.
\end{enumerate}
\end{lemma}

\begin{lemma} \label{lem:case2}
Let
$\Gamma$ denote a distance-regular graph with diameter
$D\ge 3$ and valency $k$.
Assume
 $a_i = 0$ for $0 \leq i \leq D-2$ and $a_{D-1}\not=0$.
 Let $\theta$ denote any eigenvalue of the matrix (\ref{eq:mat}).
Then the following (i)--(iv) hold.
\begin{enumerate}
\item
$\theta \in \R$.
\item Let
 $\sigma_0, \sigma_1, \ldots, \sigma_D$
denote the  pseudo cosine sequence for $\theta$. Then
$\sigma_{D-1}=\sigma_D$. \item The pair $\theta, -k$ is tight.
\item The auxiliary parameter for $\theta, -k$ is $1$.
\end{enumerate}
\end{lemma}
\begin{proof}
(i) Denote the matrix (\ref{eq:mat}) by $G$.  Observe $G$ is tridiagonal. The entries on the superdiagonal
and subdiagonal of $G$ are positive.
Therefore there exists a diagonal matrix $N$ such that
$NGN^{-1}$  is symmetric and has all entries real.
We  now see the eigenvalues of
$NGN^{-1}$ are real.
The matrices $G$ and
$NGN^{-1}$ have the same eigenvalues, so the eigenvalues
of $G$ are real.
\\
\noindent (ii) The requirement that
$(\sigma_0, \sigma_1, \ldots, \sigma_{D-1})^t$
be an eigenvector of the matrix (\ref{eq:mat}) with $\sigma_0=1$ is  the same
as the requirement that
 $\sigma_0, \sigma_1, \ldots, \sigma_{D-1},\sigma_{D-1}$
 be the pseudo cosine sequence for $\theta$.
\\
\noindent (iii),(iv) 
Let
 $\sigma_0, \sigma_1, \ldots, \sigma_D$
denote the  pseudo cosine sequence for $\theta$ and
let $\rho_0, \rho_1, \ldots, \rho_D$ denote the pseudo cosine sequence
for $-k$. 
Define $\varepsilon = 1$.
By Lemma
\ref{lem:pcformk}, we have $\rho_i = (-1)^i$ for $0 \leq i \leq D-1$.
By this and since $\sigma_{D-1}=\sigma_D$
we find
(\ref{eq:epseq}) holds. Applying
Theorem
\ref{thm:main1}
we find
 $\sigma_0, \sigma_1, \ldots, \sigma_D$
and
$\rho_0, \rho_1, \ldots, \rho_D$ form a tight pair.
Now the pair $\theta, -k$ is tight. By construction the auxiliary parameter for this pair
is $1$.
\end{proof}

\begin{theorem}\label{thm:case2} 
Let
$\Gamma$ denote a distance-regular graph with diameter
$D\ge 3$ and valency $k$.
Assume
 $a_i = 0$ for $0 \leq i \leq D-2$ and $a_{D-1}\not=0$.
Then the following (i)--(iii) hold.
\begin{enumerate}
\item The pair $\theta, k$ is tight for all $\theta \in \R$.
\item The pair $\theta, -k$ is tight for any eigenvalue $\theta$ of the matrix (\ref{eq:mat}). 
\item $\Gamma$ has no further tight pairs.
\end{enumerate}
\end{theorem}
\begin{proof}(i)  This is from Example
\ref{ex:triv}.
\\
\noindent (ii) Follows from Lemma \ref{lem:case2}(ii).
\\
\noindent (iii) Let $\theta, \theta'$ denote a tight pair,
such that $\theta \not=k$, $\theta'\not=k$.
Applying Lemma
\ref{lem:tpev}
we find $\theta, \theta'$ satisfy
(\ref{eq:ththp}). By this and since  $a_1=0$
we find $(\theta+k)(\theta'+k)=0$.
Now $\theta=-k$ or
$\theta'=-k$.
Without loss of generality we assume
$\theta'=-k$.
Let $\varepsilon$ denote the auxiliary parameter for $\theta, -k$.
By
(\ref{eq:auxform}) (with $\rho=-1$) we find $\varepsilon=1$.
By Lemma \ref{lem:Gmat}
and Lemma
\ref{lem:threeres}(ii) we find
$\theta$ is an eigenvalue of the matrix (\ref{eq:mat}).
The result follows.
\end{proof}

\section{The case $a_1=0$ and there exists $i$ $(2 \leq i \leq D-2)$
such that $a_i\not=0$}

\begin{theorem}\label{thm:case3}
Let
$\Gamma$ denote a distance-regular graph with diameter
$D\ge 4$ and valency $k$.
Assume
 $a_1 = 0$ and
there exists $i$ $(2 \leq i \leq D-2)$ such that $a_i\not=0$.
\begin{enumerate}
\item The pair $\theta, k$ is tight for all $\theta \in \R$.
\item $\Gamma$ has no further tight pairs.
\end{enumerate}
\end{theorem}

\begin{proof}
Let $\theta, \theta'$ denote a tight pair,
such that $\theta \not=k$, $\theta'\not=k$.
Applying Lemma
\ref{lem:tpev}
we find $\theta, \theta'$ satisfy
(\ref{eq:ththp}). By this and since  $a_1=0$
we find $(\theta+k)(\theta'+k)=0$.
Now $\theta=-k$ or
$\theta'=-k$.
Without loss of generality we assume
$\theta'=-k$.
Let $\varepsilon$ denote the auxiliary parameter for $\theta, -k$.
By
(\ref{eq:auxform}) (with $\rho=-1$) we find $\varepsilon=1$.
By Lemma \ref{lem:threeres}(iii) we find
$a_i=0$ for $0 \leq i \leq D-2$, a contradiction.
\end{proof}

\section{The case $a_1\not=0$}

Let $\Gamma$ denote a distance-regular graph with diameter $D\ge
3$ and $a_1\not=0$. Suppose we are given nontrivial pseudo cosine
sequences
 $\sigma_0, \sigma_1, \ldots, \sigma_D$
and $\rho_0, \rho_1, \ldots, \rho_D$ which form a tight pair. Let
$\varepsilon $ denote the corresponding auxiliary parameter. Pick an integer $i\
(1 \le i \le D)$.  Note (\ref{eq:epseq}) holds; rearranging the terms in 
(\ref{eq:epseq}) we obtain
\begin{equation}\label{eq:rhoicoeff}
\rho_{i-1}(\sigma_{i-1}-\varepsilon\sigma_i)=\rho_i
(\sigma_i-\varepsilon\sigma_{i-1}) .
\end{equation}
 We would like to
solve for $\rho_{i-1}$ in (\ref{eq:rhoicoeff}). However the
coefficient $\sigma_{i-1}-\varepsilon\sigma_{i}$ might be zero. In
the following lemma we investigate this possibility.

\begin{lemma}\label{lem:zerorhoi}
Let $\Gamma$ denote a distance-regular graph with diameter $D\ge
3$ and $a_1\not=0$. Suppose we are given nontrivial pseudo cosine
sequences
 $\sigma_0, \sigma_1, \ldots, \sigma_D$
and $\rho_0, \rho_1, \ldots, \rho_D$ which form a tight pair. Let
$\varepsilon $ denote the corresponding auxiliary parameter. For
$1 \le i \le D-1$ the following (i)--(iv) are equivalent:
\begin{enumerate}
\item $\sigma_{i-1}=\varepsilon \sigma_i.$ \item
$\sigma_{i+1}=\varepsilon \sigma_i.$ \item
$\sigma_{i-1}=\sigma_{i+1}.$ \item $\rho_i =0.$
\end{enumerate}
\end{lemma}

\begin{proof} Recall line
(\ref{eq:epseq}) holds.\\
(i)$\Longrightarrow$ (iv) We replace $\sigma_{i-1}$ by
$\varepsilon\sigma_i$ in (\ref{eq:epseq}) to obtain
$\sigma_i\rho_i(1-\varepsilon^2)=0.$ Observe
 $\varepsilon  \not\in \lbrace 1, -1 \rbrace$
by Corollary \ref{cor:epsgood}. Assume for the moment that
$\sigma_i=0$. It follows from the assumption that
$\sigma_{i-1}=0.$  So $\sigma_{i-1}=\sigma_i$ contradicting
Corollary \ref{cor:epsgood}.
 Thus, $\sigma_i\ne 0$ so $\rho_i=0$.\\
(iv)$\Longrightarrow$ (i) Setting  $\rho_i=0$ in (\ref{eq:epseq})
we obtain $\rho_{i-1}(\sigma_{i-1} -\varepsilon\sigma_i)=0.$
Observe $\rho_{i-1}\ne 0$ by Lemma \ref{lem:doublezero} so
$\sigma_{i-1}=\varepsilon \sigma_i.$\\
(ii)$\Longleftrightarrow$(iv)  Similar to the proof of (i)$\Longleftrightarrow$(iv).\\
(i),(ii)$\Longrightarrow$(iii)  Clear.\\
(iii)$\Longrightarrow$(i) Adding (\ref{eq:epseq}) at $i$ and $i+1$
we obtain
$$
\sigma_{i+1}\rho_{i+1}- \sigma_{i-1}\rho_{i-1}=
\varepsilon(\sigma_i\rho_{i+1} -\sigma_{i+1}\rho_i
+\sigma_{i-1}\rho_i -\sigma_i\rho_{i-1}).
$$
Replacing $\sigma_{i+1}$ by $\sigma_{i-1}$ in the above line and
simplifying, we obtain
\begin{equation}\label{eq:tech}
(\sigma_{i-1}-\varepsilon \sigma_i)(\rho_{i+1}-\rho_{i-1})=0.
\end{equation}
We show $\rho_{i+1}-\rho_{i-1}\neq 0$. Suppose
$\rho_{i+1}=\rho_{i-1}$. Using this and Lemma \ref{lem:triplemult}
we have $\rho_{i-1}-\rho_i\ne 0.$  Applying  Lemma
\ref{lem:recur2} to  $\rho_0, \rho_1, \ldots, \rho_D$ we obtain
\begin{equation}\label{eq:aik1}
\frac{a_i}{k}=\frac{\rho_{i-1}-\rho\rho_i}{\rho_{i-1}-\rho_i}.
\end{equation}
Similarly since $\sigma_{i+1}=\sigma_{i-1}$ we have
$\sigma_{i-1}-\sigma_i\ne 0$ and
\begin{equation}\label{eq:aik2}
\frac{a_i}{k}=\frac{\sigma_{i-1}-\sigma\sigma_i}{\sigma_{i-1}-\sigma_i}.
\end{equation}
Combining (\ref{eq:aik1}) and (\ref{eq:aik2}) we obtain
\begin{equation}\label{eq:aieq}
(\sigma_{i-1}-\sigma\sigma_i)(\rho_{i-1}-\rho_i)=
(\sigma_{i-1}-\sigma_i)(\rho_{i-1}-\rho\rho_i).
\end{equation}
We view (\ref{eq:epseq}) and (\ref{eq:aieq})  as a homogeneous
system of linear equations in $\sigma_{i-1}$ and $\sigma_i$. We
find that the coefficient matrix of this linear system is given by
$$C=\left(\begin{array}{cc}
\varepsilon\rho_i +\rho_{i-1} & -\rho_i-\varepsilon\rho_{i-1}\\
(\rho-1)\rho_i & \rho_{i-1}(1-\sigma)+\rho_i(\sigma-\rho)
\end{array}\right).$$
Observe $\sigma\ne \rho$ by Lemma \ref{lem:dist} so
(\ref{eq:auxform}) holds.  Computing the determinant of $C$ and
evaluating the result using (\ref{eq:auxform}) we find ${\rm det}\
C=(\rho_{i-1}-\rho_i)(1-\sigma)(\rho_{i-1}-\rho\rho_i).$ We
mentioned earlier that $\rho_{i-1}-\rho_i \not= 0$.  Observe $1-
\sigma \ne 0$ since the sequence $\sigma_0, \sigma_1,\ldots,
\sigma_D$ is nontrivial. Observe $\rho_{i-1}-\rho\rho_i\ne 0$;
otherwise $a_i=0$ in view of (\ref{eq:aik1}), contradicting Lemma
\ref{lem:a1notzero}. We now see that ${\rm det}\ C\not=0$ so $C$
is nonsingular.  Therefore  $\sigma_{i-1}=0$ and $\sigma_i=0$.
This contradicts Lemma \ref{lem:doublezero}. Thus we have shown
$\rho_{i+1}-\rho_{i-1}\neq 0$. By this and (\ref{eq:tech}) we find
$\sigma_{i-1}=\varepsilon \sigma_i$.
\end{proof}

\begin{corollary}\label{cor:rhominusplus}
Let $\Gamma$ denote a distance-regular graph with diameter $D\ge
3$ and $a_1\not=0$. Suppose we are given nontrivial pseudo cosine
sequences
 $\sigma_0, \sigma_1, \ldots, \sigma_D$
and $\rho_0, \rho_1, \ldots, \rho_D$ which form a tight pair. Let
$\varepsilon $ denote the corresponding auxiliary parameter. Pick
an integer $i\ (1 \le i \le D-1)$  and assume  the equivalent
conditions (i)--(iv) of Lemma \ref{lem:zerorhoi} do not hold. Then
both
\begin{eqnarray}
\rho_{i-1}&=&\rho_i\frac{\sigma_i-\varepsilon\sigma_{i-1}}{\sigma_{i-1}-\varepsilon\sigma_i},\label{eq:rhominus}\\
\rho_{i+1} &=&
\rho_i\frac{\sigma_i-\varepsilon\sigma_{i+1}}{\sigma_{i+1}-\varepsilon\sigma_i}.\label{eq:rhoplus}
\end{eqnarray}
\end{corollary}

\begin{proof}
 Observe (\ref{eq:epseq})  holds.   Rearranging terms we find
\begin{equation}\label{eq:faveqati}
\rho_{i-1}(\sigma_{i-1}-\varepsilon\sigma_i)=\rho_i(\sigma_i-\varepsilon\sigma_{i-1}).
\end{equation}
Observe $\sigma_{i-1}-\varepsilon\sigma_i\not=0$ by Lemma
\ref{lem:zerorhoi}(i); solving (\ref{eq:faveqati})
for $\rho_{i-1}$ we obtain (\ref{eq:rhominus}). 
Replacing $i$ by $i+1$ in 
(\ref{eq:faveqati}) and rearranging terms we obtain
\begin{equation}\label{eq:faveqatip}
\rho_{i+1}(\sigma_{i+1}-\varepsilon\sigma_{i})=\rho_{i}(\sigma_{i}-\varepsilon\sigma_{i+1}).
\end{equation}
Observe  $\sigma_{i+1}-\varepsilon\sigma_{i}\not=0$ by Lemma
\ref{lem:zerorhoi}(ii).  Solving  (\ref{eq:faveqatip}) for
$\rho_{i+1}$ we get (\ref{eq:rhoplus}).
\end{proof}

\section{Tight pseudo cosine sequences}
Let $\Gamma$ denote a distance-regular graph with diameter $D\ge
3$ and $a_1\ne 0$.  Let  $\sigma_0, \sigma_1, \ldots, \sigma_D$
denote a nontrivial pseudo cosine sequence.  We want to prove there exists at most one nontrivial
pseudo cosine sequence $\rho_0, \rho_1, \ldots, \rho_D$ such that $\sigma_0,
\sigma_1, \dots, \sigma_D$ and $\rho_0, \rho_1, \ldots, \rho_D$
 form a tight pair.
To do this we need a lemma.

\begin{lemma}\label{lem:unique}
Let $\Gamma$ denote a distance-regular graph with diameter $D\ge
3$ and valency $k$. Suppose there exist real numbers $\theta, \theta', \theta''$ such that (i) none of
$\theta, \theta', \theta''$  is equal to $k$; (ii)
$\theta' \ne \theta''$; (iii) the pair $\theta, \theta'$ is tight; and
(iv) the pair $\theta, \theta''$ is tight. Then
$a_i=0$ for $0 \le i \le D-2.$  Moreover $\theta=-k$.
\end{lemma}

\begin{proof}
By Lemma \ref{lem:tpev} the pair $\theta, \theta'$  satisfies
(\ref{eq:ththp}). Applying the same lemma to the pair $\theta,
\theta''$ and using the fact that $\theta' \ne \theta''$ we
obtain $\theta+k/(a_1+1)=0$. Now in (\ref{eq:ththp}) the left-hand side
 is $0$ so the right-hand side is $0$.  Therefore $a_1=0$ and $\theta =-k$.
 If $D=3$, we are done.  Assume $D\ge 4$.  By  Theorem \ref{thm:case3} we have
$a_2=\cdots = a_{D-2}=0$. The result follows.
\end{proof}

\begin{corollary}
Let $\Gamma$ denote a distance-regular graph with diameter $D\ge
3$ and $a_1\ne 0$.
Let $\sigma_0, \sigma_1, \dots, \sigma_D$
denote a nontrivial pseudo cosine sequence. Then there exists
at most one nontrivial pseudo cosine sequence
$\rho_0, \rho_1, \ldots, \rho_D$ such that $\sigma_0,
\sigma_1, \dots, \sigma_D$ and $\rho_0, \rho_1, \ldots, \rho_D$
 form a tight pair. Suppose $\rho_0, \rho_1, \ldots, \rho_D$ exists. 
 Then the corresponding auxiliary parameter is unique.
\end{corollary}

\begin{proof}
The first assertion follows from
Lemma \ref{lem:unique}. The auxiliary parameter 
is unique by Lemma \ref{lem:dist} and the comment at the end of
Definition \ref{def:auxpar}.
\end{proof}

\begin{definition}\label{def:tightpseudo}
{\rm 
Let $\Gamma$ denote a distance-regular graph with diameter $D\ge
3$ and $a_1\ne 0$. Let $\sigma_0, \sigma_1, \dots, \sigma_D$
denote a nontrivial pseudo cosine sequence.  We say this sequence
is \emph{tight} whenever there exists a nontrivial pseudo cosine
sequence $\rho_0, \rho_1, \ldots, \rho_D$ such that $\sigma_0,
\sigma_1, \dots, \sigma_D$ and $\rho_0, \rho_1, \ldots, \rho_D$
 form a tight pair.  By the \emph{auxiliary parameter} for
 $\sigma_0, \sigma_1, \dots, \sigma_D$ we mean the auxiliary
 parameter for the tight pair $\sigma_0,
\sigma_1, \dots, \sigma_D$ and $\rho_0, \rho_1, \ldots, \rho_D$.}
\end{definition}

\begin{lemma}\label{lem:fillgap}
Let $\Gamma$ denote a distance-regular graph with diameter $D\ge
3$ and $a_1\not=0$. Suppose we are given a nontrivial pseudo
cosine sequence $\sigma_0, \sigma_1, \ldots, \sigma_D$ which is
tight in the sense of Definition \ref{def:tightpseudo}. Let $\varepsilon $ denote the corresponding auxiliary parameter.  Then the following (i)--(iv) hold.
\begin{enumerate}
\item $\sigma_2 \ne 1$. 
\item $\varepsilon\sigma\not= 1$. 
\item $\sigma_2 \not= \varepsilon\sigma$. 
\item $\sigma_2 \ne \sigma^2$.
\end{enumerate}
\end{lemma}

\begin{proof}
By Definition \ref{def:tightpseudo} there exists a nontrivial
pseudo cosine sequence $\rho_0, \rho_1, \ldots, \rho_D$ such that
$\sigma_0, \sigma_1, \dots, \sigma_D$ and $\rho_0, \rho_1, \ldots,
\rho_D$ form a tight pair. Thus (\ref{eq:epseq}) holds.\\
(i) Setting $i=3$ in (\ref{eq:epseq}) we find
\begin{equation}\label{eq:epseqat3}
\sigma_3\rho_3-\sigma_2\rho_2-\varepsilon(\sigma_2\rho_3-\sigma_3\rho_2)
\end{equation} is zero. We assume $\sigma_2=1$ and show (\ref{eq:epseqat3}) is
not zero.  In order to do this we evaluate the terms in
(\ref{eq:epseqat3}). We assume $\sigma_2=1$; therefore the equivalent conditions
(i)--(iv) hold in Lemma \ref{lem:zerorhoi} for $i=1$.  Setting $i=1$ in Lemma \ref{lem:zerorhoi}(iv) we find $\rho=0$.  Similarly using Lemma \ref{lem:zerorhoi}(i)
we find $\varepsilon\sigma=1.$
It follows that $\sigma\not= 0$ and $\varepsilon=1/\sigma.$
Setting $i=1$ in (\ref{eq:recurc2}), solving for $\sigma$ and
eliminating $a_1$ in the result using (\ref{eq:intnos}) we get
\begin{equation}\label{eq:s1}
\sigma=-\frac{1+b_1}{k}.
\end{equation}
By this and since $\varepsilon=1/\sigma$ we find
\begin{equation}\label{eq:ep}
\varepsilon=-\frac{k}{1+b_1}.
\end{equation}
 Setting $i=2$ in (\ref{eq:recursion}), solving for $\sigma_3$ and
 eliminating $c_2,\sigma$ in the result using (\ref{eq:intnos}), (\ref{eq:s1}) we get
\begin{equation}\label{eq:s3}
\sigma_3=-\frac{(a_2+b_2)(1+b_1)+ka_2}{kb_2}.
\end{equation}
Applying (\ref{eq:recursion}) to the sequence $\rho_0, \rho_1,
\ldots, \rho_D$ we find
\begin{equation}\label{eq:rhorecursion}
c_i\rho_{i-1}+a_i\rho_i+b_i\rho_{i+1}=k\rho\rho_i \quad (0 \le i
\le D-1).
\end{equation}
 Setting $i=1,2$ in (\ref{eq:rhorecursion}) and
solving for $\rho_2, \rho_3$ respectively using $\rho=0$ we obtain
\begin{eqnarray}
\rho_2 &=& -{b_1}^{-1},\label{eq:r2}\\
\rho_3 &=& {a_2}{b_1}^{-1}{b_2}^{-1}.\label{eq:r3}
\end{eqnarray}
Evaluating (\ref{eq:epseqat3}) using  (\ref{eq:s1}), (\ref{eq:ep}), (\ref{eq:s3}) and (\ref{eq:r2}), (\ref{eq:r3}) we routinely verify that
 \begin{equation}\label{eq:contradiction}
\sigma_3\rho_3-\sigma_2\rho_2-\varepsilon
(\sigma_2\rho_3-\sigma_3\rho_2)=-a_2\frac{(a_2+b_2)(1+b_1+k)}{b_1b_2^2k}.
\end{equation}
Recall the left-hand side of (\ref{eq:contradiction}) is $0$ so the right-hand side is $0.$ 
It follows $a_2=0$.  However $a_2\ne 0$ by Lemma \ref{lem:a1notzero} and since $a_1 \not=0$. We have now shown $\sigma_2\ne 1$.\\
 (ii),(iii) Immediate from (i) above and Lemma \ref{lem:zerorhoi}.\\
 (iv) Suppose $\sigma_2=\sigma^2$.  Setting $i=1$ in
 (\ref{eq:recurc2}) and since $\sigma\not= 1$ we find $\sigma=-1/(a_1+1)$. Setting
 $\theta=k\sigma, \theta'=k\rho$ in equation (\ref{eq:ththp}) we find that the
 left-hand side is $0$ while the right-hand side is nonzero. This
 violates Lemma \ref{lem:tpev}.  Thus $\sigma_2\not=\sigma^2.$
\end{proof}

\begin{lemma}\label{lem:rhorho2}
Let $\Gamma$ denote a distance-regular graph with diameter $D\ge
3$ and $a_1\not=0$. Suppose we are given nontrivial pseudo cosine
sequences
 $\sigma_0, \sigma_1, \ldots, \sigma_D$
and $\rho_0, \rho_1, \ldots, \rho_D$ which form a tight pair. Let
$\varepsilon $ denote the corresponding auxiliary parameter. Then
the following (i)--(iii) hold.
\begin{enumerate}
\item $\sigma\not= \varepsilon.$ 
\item $\rho = \displaystyle{\frac{1-\varepsilon\sigma}{\sigma-\varepsilon}}.$
\item $\rho_2 =\displaystyle{\frac{\rho(\sigma-\varepsilon\sigma_2)}{\sigma_2-\varepsilon\sigma}}.$
\end{enumerate}
\end{lemma}

\begin{proof} (i) Suppose $\sigma = \varepsilon$.  Then setting
 $i=1$ in (\ref{eq:epseq}) we find
$\varepsilon^2 = 1;$ contradicting Corollary \ref{cor:epsgood}.\\
(ii) Solve (\ref{eq:auxform}) for $\rho$.\\
(iii) Set $i=1$ in (\ref{eq:rhoplus}).
\end{proof}

\section{The intersection numbers}

Let $\Gamma$ denote a distance-regular graph with diameter $D\ge
3$ and $a_1\not=0$.  Suppose we are given a nontrivial pseudo
cosine sequence
 $\sigma_0, \sigma_1, \ldots, \sigma_D$
which is tight. Let $\varepsilon $ denote the corresponding
auxiliary parameter.  In this section we compute the intersection numbers
in terms of  $\sigma_0, \sigma_1, \ldots, \sigma_D$ and $\varepsilon.$ \\

\noindent We begin with the valency $k$.

\begin{lemma}\label{lem:kform}
 Let $\Gamma$ denote a distance-regular graph with diameter $D\ge
3$ and $a_1\not=0$. Suppose we are given a nontrivial pseudo
cosine sequence
 $\sigma_0, \sigma_1, \ldots, \sigma_D$
which is tight. Let $\varepsilon $ denote the corresponding
auxiliary parameter.  Then the valency $k$ satisfies
\begin{equation}\label{eq:kformA}
k = h\frac{\sigma-\varepsilon}{\sigma - 1}
\end{equation}
where
\begin{equation}\label{eq:hformA}
h = \frac{(1-\sigma)(1-\sigma_2)}{(\sigma^2-\sigma_2)
(1-\varepsilon \sigma)}.
\end{equation}
We remark the denominator in (\ref{eq:kformA}) is nonzero since
$\sigma_0, \sigma_1, \ldots, \sigma_D$ is nontrivial.  Moreover the
denominator in (\ref{eq:hformA}) is nonzero by Lemma
\ref{lem:fillgap}.
\end{lemma}

\begin{proof}
By Definition \ref{def:tightpseudo} there exists a nontrivial pseudo cosine sequence
$\rho_0, \rho_1, \ldots, \rho_D$ such that $\sigma_0, \sigma_1,
\dots, \sigma_D$ and $\rho_0, \rho_1, \ldots, \rho_D$ form a tight
pair.  Since $\sigma_0, \sigma_1, \dots, \sigma_D$ is a pseudo
cosine sequence we find it satisfies (\ref{eq:recur}). Setting
$i=1$ in  (\ref{eq:recur})  we get
 \begin{equation}
k(\sigma^2-\sigma)+b_1(\sigma-\sigma_2) =
1-\sigma.\label{eq:sigrecur}
\end{equation}
Similarly\
\begin{equation}
 k(\rho^2-\rho)+b_1(\rho-\rho_2) =
1-\rho.\label{eq:rhorecur}
 \end{equation}
We view (\ref{eq:sigrecur}) and (\ref{eq:rhorecur})  as a system
of linear equations in $k$ and $b_1$. The coefficient matrix is
$$E=\left(\begin{array}{cc}
 \sigma^2-\sigma & \sigma-\sigma_2 \\
 \rho^2-\rho & \rho-\rho_2
\end{array}\right).$$
 Evaluating the determinant of $E$ using Lemma \ref{lem:rhorho2}(ii),(iii) we find that
$${\rm
det}\ E =
\frac{(1+\varepsilon)(1-\sigma)(1-\varepsilon\sigma)(\sigma-\sigma_2)(\sigma^2-\sigma_2)}
{(\sigma-\varepsilon)^2(\sigma_2-\varepsilon\sigma)}.$$ 
We show ${\rm
det}\ E \not=0$.  Observe $1+\varepsilon\not= 0$ by Corollary \ref{cor:epsgood},
$1-\sigma\ne 0$ since the sequence $\sigma_0,
\sigma_1, \dots, \sigma_D$ is nontrivial, $1-\varepsilon\sigma\not=
0$ by Lemma \ref{lem:fillgap}(ii),  $\sigma-\sigma_2\ne 0$ by Corollary \ref{cor:epsgood} and $\sigma^2-\sigma_2\ne 0$
by Lemma \ref{lem:fillgap}(iv). We have now shown ${\rm det}\ E\not= 0$. 
 Solving the system (\ref{eq:sigrecur}), (\ref{eq:rhorecur})  we find 
$$k=-\frac{(\sigma-\varepsilon)(1-\sigma_2)}{(\sigma^2-\sigma_2)(1-\varepsilon\sigma)}.
$$
The desired result follows.
\end{proof}

\begin{lemma}\label{lem:aiform}
 Let $\Gamma$ denote a distance-regular graph with diameter $D\ge
3$ and $a_1\not=0$. Suppose we are given a nontrivial pseudo
cosine sequence
 $\sigma_0, \sigma_1, \ldots, \sigma_D$
which is tight. Let $\varepsilon $ denote the corresponding
auxiliary parameter.  Then
\begin{equation}\label{eq:aiformA}
a_i = g \frac{(\sigma_{i+1}-\sigma \sigma_i) (\sigma_{i-1}-\sigma
\sigma_i)}{ (\sigma_{i+1}- \sigma_i) (\sigma_{i-1}- \sigma_i)}
\qquad \qquad (1 \leq i \leq D-1)
\end{equation}
where
\begin{eqnarray}\label{eq:gformA}
g = \frac{(\varepsilon-1)(1-\sigma_2)}
{(\sigma^2-\sigma_2)(1-\varepsilon \sigma)}.
\end{eqnarray}
We remark the denominator in (\ref{eq:aiformA}) is nonzero by
Corollary \ref{cor:epsgood}. Moreover the denominator in
(\ref{eq:gformA}) is nonzero by Lemma \ref{lem:fillgap}.
\end{lemma}

\begin{proof}
By Definition \ref{def:tightpseudo}  there exists a nontrivial pseudo cosine sequence
$\rho_0, \rho_1, \ldots, \rho_D$ such that $\sigma_0, \sigma_1,
\dots, \sigma_D$ and $\rho_0, \rho_1, \ldots, \rho_D$ form a tight
pair.  Pick an integer $i\ (1\le i \le D-1).$
By Lemma \ref{lem:recur2} and since $\sigma_0, \sigma_1,
\dots, \sigma_D$ is a pseudo cosine sequence we find
\begin{equation}\label{eq:recursigma}
c_i(\sigma_{i+1}-\sigma_{i-1})+a_i(\sigma_{i+1}-\sigma_i) =
k(\sigma_{i+1}-\sigma\sigma_i).
\end{equation}
Similarly
\begin{equation}\label{eq:recurrho}
c_i(\rho_{i+1}-\rho_{i-1})+a_i(\rho_{i+1}-\rho_i) =
k(\rho_{i+1}-\rho\rho_i).
\end{equation}
We solve (\ref{eq:recursigma}), (\ref{eq:recurrho}) for $a_i$.
 We consider two cases.  First assume $\sigma_{i+1}-\sigma_{i-1}\ne 0$.  We view (\ref{eq:recursigma})
and (\ref{eq:recurrho}) as a system of linear equations in $c_i$
and $a_i$.  The coefficient matrix of the above system is given by
$$F=\left(\begin{array}{cc}
 \sigma_{i+1}-\sigma_{i-1} & \sigma_{i+1}-\sigma_i\\
 \rho_{i+1}-\rho_{i-1} & \rho_{i+1}-\rho_i
\end{array}\right).$$
We evaluate the determinant of $F$. We assume $\sigma_{i+1}\ne\sigma_{i-1}$ so
 Corollary \ref{cor:rhominusplus} applies; therefore
$\rho_{i-1}, \rho_{i+1}$ in $F$ are given by (\ref{eq:rhominus}),
(\ref{eq:rhoplus}) respectively. 
 Eliminating $\rho_{i-1}$ and $\rho_{i+1}$
in $F$ using (\ref{eq:rhominus}) and (\ref{eq:rhoplus}) we get
$${\rm det}\
F=\rho_i\frac{(1+\varepsilon)(\sigma_{i-1}-\sigma_i)(\sigma_{i+1}-\sigma_i)(\sigma_{i-1}-\sigma_{i+1})}
{(\sigma_{i-1}-\varepsilon\sigma_i)(\sigma_{i+1}-\varepsilon\sigma_i)}.
$$
 We show ${\rm det}\ F\not=0.$  Observe $\rho_i \not = 0$ by Lemma \ref{lem:zerorhoi}
 and since $\sigma_{i+1}\ne\sigma_{i-1}$. Observe $1+\varepsilon\not= 0$,
 $\sigma_{i-1}-\sigma_i\not= 0$ and $\sigma_{i+1}-\sigma_i\not= 0$ by Corollary \ref{cor:epsgood}.  Now ${\rm det}\ F\not=0.$
 Solving (\ref{eq:recursigma}), (\ref{eq:recurrho})
 for $a_i$ and then evaluating the result using Corollary \ref{cor:rhominusplus},
 Lemma \ref{lem:rhorho2}(ii) and Lemma \ref{lem:kform}
 we routinely obtain (\ref{eq:aiformA}), (\ref{eq:gformA}).
 We have now proved the result for the case $\sigma_{i+1}-\sigma_{i-1}\not=0$. Next assume
 $\sigma_{i+1}-\sigma_{i-1}=0$. Observe
 $\sigma_{i+1}-\sigma_i\ne 0$ by  Lemma \ref{lem:triplemult}.  Solving for $a_i$ in
  (\ref{eq:recursigma})  we get
 \begin{equation}\label{eq:case1}\\
a_i=k\frac{\sigma_{i+1}-\sigma\sigma_i}{\sigma_{i+1}-\sigma_i}.
 \end{equation}
Eliminating $k$ in (\ref{eq:case1}) using (\ref{eq:kformA}) and
(\ref{eq:hformA}) we get
 \begin{equation}\label{eq:case1b}\\
a_i=\frac{(1-\sigma_2)(\varepsilon-\sigma)(\sigma_{i+1}-\sigma\sigma_i)}{(\sigma^2-\sigma_2)(1-\varepsilon\sigma)(\sigma_{i+1}-\sigma_i)}.
 \end{equation}
 We eliminate $\varepsilon -\sigma$ in (\ref{eq:case1b}).
By Lemma
 \ref{lem:zerorhoi} and since  $\sigma_{i+1}=\sigma_{i-1}$ we find  $\sigma_{i-1}=\varepsilon\sigma_i$.  Observe
$\sigma_{i-1}-\sigma_i \not= 0$ by 
 Corollary \ref{cor:epsgood}. Using these facts we find
 \begin{equation}\label{eq:link}
\varepsilon-\sigma =\frac{(\varepsilon-1)
(\sigma_{i-1}-\sigma\sigma_i)}{\sigma_{i-1}-\sigma_i}.
\end{equation}
Eliminating $\varepsilon-\sigma$ in (\ref{eq:case1b}) using
(\ref{eq:link}) we get (\ref{eq:aiformA}),
(\ref{eq:gformA}).
\end{proof}

\begin{lemma}\label{lem:kbiciform}
 Let $\Gamma$ denote a distance-regular graph with diameter $D\ge
3$ and $a_1\not=0$. Suppose we are given a nontrivial pseudo
cosine sequence
 $\sigma_0, \sigma_1, \ldots, \sigma_D$
which is tight. Let $\varepsilon $ denote the corresponding
auxiliary parameter.  Then
\begin{eqnarray}
b_i (\sigma_{i-1}- \sigma_{i+1}) &=& h \frac{(\sigma_{i-1}-\sigma
\sigma_i) (\sigma_{i+1}-\varepsilon \sigma_i)} {\sigma_{i+1}-
\sigma_{i}} \qquad (1 \leq i \leq D-1), \label{eq:biformA}\\
c_i (\sigma_{i+1}- \sigma_{i-1}) &=& h \frac{(\sigma_{i+1}-\sigma
\sigma_i) (\sigma_{i-1}-\varepsilon \sigma_i)} {\sigma_{i-1}-
\sigma_{i}} \qquad (1 \leq i \leq D-1), \label{eq:ciformA}
\end{eqnarray}
where $h$ is as in (\ref{eq:hformA}). We remark  
the denominators in (\ref{eq:biformA}), (\ref{eq:ciformA}) are
nonzero by Corollary \ref{cor:epsgood}.
\end{lemma}

\begin{proof}
To obtain (\ref{eq:biformA}), in line (\ref{eq:recurb2}) evaluate $k, a_i$ using
Lemmas \ref{lem:kform},
\ref{lem:aiform} respectively.  Line (\ref{eq:ciformA}) is similarly obtained.
\end{proof}

\section{A characterization theorem}

\begin{definition}\label{def:condition}
{\rm
Let $\Gamma$ denote a distance-regular graph with diameter $D\ge
3$ and $a_1 \ne 0$. Let $\sigma_0, \sigma_1,\ldots, \sigma_D$ and
$\varepsilon$ denote real numbers such that $\sigma_0=1$.  We
consider several conditions on these scalars. 
\begin{enumerate}
\item By condition $A$ we
mean $\varepsilon\not= -1$, equations (\ref{eq:aiformA}),
(\ref{eq:gformA}) hold and the denominators in
(\ref{eq:aiformA}), (\ref{eq:gformA}) are nonzero. 
\item By condition
$B$ we mean equations (\ref{eq:biformA}), (\ref{eq:hformA}) hold
and the denominators in (\ref{eq:biformA}),
(\ref{eq:hformA}) are nonzero. 
\item By condition $C$ we mean equations
(\ref{eq:kformA}), (\ref{eq:ciformA}), (\ref{eq:hformA}) hold and
the denominators in (\ref{eq:kformA}),
(\ref{eq:ciformA}),(\ref{eq:hformA}) are nonzero.
\end{enumerate}}
\end{definition}

\begin{theorem}
\label{th:mainth}
Let
$\Gamma$ denote a distance-regular graph with diameter
$D\ge 3$ and $a_1 \not=0$.
Let $\sigma_0, \sigma_1,\ldots, \sigma_D$ and $\varepsilon$
denote real numbers such that $\sigma_0=1$.
Then with reference to Definition \ref{def:condition} the following (i)--(iv)
are equivalent.
\begin{enumerate}
\item
 $\sigma_0, \sigma_1,\ldots, \sigma_D$ is a tight nontrivial pseudo
cosine sequence
 and $\varepsilon$ is the corresponding auxiliary parameter.
\item
 $\sigma_0, \sigma_1,\ldots, \sigma_D$ is a  nontrivial pseudo
cosine sequence and
 $\sigma_0, \sigma_1,\ldots, \sigma_D, \varepsilon$
satisfy $A$.
 \item
 $\sigma_0, \sigma_1,\ldots, \sigma_D, \varepsilon$ satisfy
 both $A$ and $B$.
\item
 $\sigma_0, \sigma_1,\ldots, \sigma_D, \varepsilon$ satisfy
both $A$ and $C$.
\end{enumerate}
\end{theorem}

\begin{proof}
(i)$\Longrightarrow$(ii) Applying Definition \ref{def:tightpseudo} there exists a nontrivial pseudo cosine sequence $\rho_0, \rho_1, \ldots, \rho_D$ such that
$\sigma_0,
\sigma_1,\ldots, \sigma_D$ and
 $\rho_0, \rho_1,\ldots, \rho_D$ form a tight pair.
Applying Corollary \ref{cor:epsgood} to this pair we find $\varepsilon \ne -1$.
The result follows by Lemma \ref{lem:aiform}.\\
(ii) $\Longrightarrow$ (iii) We verify that
 $\sigma_0, \sigma_1,\ldots, \sigma_D, \varepsilon$ satisfy
 condition $B$. In each of (\ref{eq:biformA}), (\ref{eq:hformA}) the denominator is
 nonzero because each factor in this denominator is in the denominator of 
 (\ref{eq:aiformA}), (\ref{eq:gformA}).  We now verify (\ref{eq:biformA}), (\ref{eq:hformA}).
 Observe Lemma
\ref{lem:recur2}(i) holds. Applying that lemma we find
(\ref{eq:recurb2}), (\ref{eq:recurc2})  hold. Setting $i=1$ in
(\ref{eq:recurc2}) we find
$k(\sigma^2-\sigma_2)=1-\sigma_2+a_1(\sigma-\sigma_2)$. Observe
$\sigma^2-\sigma_2$  is not zero since it is a factor in the
denominator of (\ref{eq:gformA}). Solving the equation for $k$ and
simplifying the result using (\ref{eq:aiformA}) we obtain
(\ref{eq:kformA}),
(\ref{eq:hformA}). Evaluating (\ref{eq:recurb2}) using this and
(\ref{eq:aiformA}) we obtain
 (\ref{eq:biformA}).
The result follows.\\
(iii) $\Longrightarrow$ (iv) We verify that
 $\sigma_0, \sigma_1,\ldots, \sigma_D, \varepsilon$ satisfy
 condition $C$. The denominator in (\ref{eq:hformA}) is nonzero by condition $B$.
  In each of (\ref{eq:kformA}), (\ref{eq:ciformA}) the denominator is
 nonzero because each factor in this denominator is in the denominator of 
 (\ref{eq:hformA}), (\ref{eq:aiformA}), (\ref{eq:gformA}), (\ref{eq:biformA}).  We now verify (\ref{eq:kformA}), (\ref{eq:ciformA}), (\ref{eq:hformA}).
In order to do this 
 we first obtain (\ref{eq:kformA}). 
 Observe that in (\ref{eq:biformA})
 the coefficient of
$b_1$ is equal to $1-\sigma_2$. This coefficient is not zero; if
it is then $g=0$ in view of (\ref{eq:gformA}) which implies
$a_1=0$ for a contradiction. Now to obtain (\ref{eq:kformA}),
simplify the right-hand side of $k=b_1+a_1+1$ using
(\ref{eq:aiformA}), 
 (\ref{eq:biformA}),(\ref{eq:gformA}), (\ref{eq:hformA}).
We now have (\ref{eq:kformA}). To obtain
 (\ref{eq:ciformA}), expand the
left-hand side using (\ref{eq:intnos}) and simplify the result
using (\ref{eq:aiformA}), (\ref{eq:kformA}),
 (\ref{eq:biformA}).
\\
(iv) $\Longrightarrow$ (ii)
We  show
 $\sigma_0, \sigma_1,\ldots, \sigma_D$
is a nontrivial pseudo cosine sequence. To do this we apply Lemma
\ref{lem:recur2}(i),(iii).  Using (\ref{eq:aiformA}), (\ref{eq:kformA}),
 (\ref{eq:ciformA}) and the fact that $g=h\frac{\varepsilon-1}{1-\sigma}$
from (\ref{eq:hformA}) and (\ref{eq:gformA}) we routinely verify
(\ref{eq:recurc2}). Observe Lemma
\ref{lem:recur2}(iii) holds. Applying that lemma we find
 $\sigma_0, \sigma_1,\ldots, \sigma_D$
is a pseudo cosine sequence.  We remark $\sigma \not=1$ since $1-\sigma$
is in the denominator of (\ref{eq:aiformA}).  Therefore
 $\sigma_0, \sigma_1,\ldots, \sigma_D$ is nontrivial.\\
(ii), (iii), (iv) $\Longrightarrow$ (i) Observe $\sigma\not=1$ since 
 $\sigma_0, \sigma_1,\ldots, \sigma_D$ is nontrivial.
 Observe $\sigma \not= \varepsilon$ by
(\ref{eq:kformA}) and since $k\not=0$. Define
\begin{equation}\label{eq:rho}
\rho=\frac{1-\varepsilon \sigma}{\sigma-\varepsilon}.
\end{equation}
 Let $\rho_0, \rho_1,\ldots, \rho_D$ denote the pseudo cosine
sequence for $k\rho $. Observe $\rho\not=1;$ otherwise
$\sigma=1$ or $\varepsilon =-1$ for a contradiction.  Therefore
 $\rho_0, \rho_1,\ldots, \rho_D$
is nontrivial. We show $\sigma_0,
\sigma_1,\ldots, \sigma_D$ and
 $\rho_0, \rho_1,\ldots, \rho_D$ form a tight pair.
To do this we apply Theorem \ref{thm:main1}. Specifically we show
that for any integer $i \  (1 \leq i \leq D)$
\begin{equation}\label{eq:indhyp}
\rho_i(\sigma_i-\varepsilon\sigma_{i-1})=\rho_{i-1}(\sigma_{i-1}-\varepsilon\sigma_i).
\end{equation}
We do this by induction on $i$.  Observe that when $i=1$
(\ref{eq:indhyp}) holds in view of (\ref{eq:rho}). Now pick an
integer $i \ (1 \leq i \leq D-1)$ and assume (\ref{eq:indhyp})
holds. We show that
$\rho_{i+1}(\sigma_{i+1}-\varepsilon\sigma_i)=\rho_i(\sigma_i-\varepsilon\sigma_{i+1}).$
This will follow provided we can show
\begin{equation}\label{eq:bitimeseq}
b_i\rho_{i+1}(\sigma_{i+1}-\varepsilon\sigma_i)-b_i\rho_i(\sigma_i-\varepsilon\sigma_{i+1})=0.
\end{equation}
We evaluate the expression on the left-hand side of  (\ref{eq:bitimeseq}) in the following
manner. In that expression, we eliminate the first occurrence of $b_i$
using (\ref{eq:recursion}) and we eliminate the second occurrence
of $b_i$ using (\ref{eq:intnos}).  We conclude the left-hand side of  (\ref{eq:bitimeseq})
is equal to $x+y+z$ where
\begin{eqnarray}
x &=&
k\rho\rho_i(\sigma_{i+1}-\varepsilon\sigma_i)-k\rho_i(\sigma_i-\varepsilon\sigma_{i+1})-
a_i\rho_i(\sigma_{i+1}-\sigma_i)(1+\varepsilon),\label{eq:x}\\
y &=&-c_i\rho_{i-1}(\sigma_{i+1}-\sigma_{i-1}),\label{eq:y}\\
z &=&
c_i\rho_i(\sigma_i-\varepsilon\sigma_{i+1})-c_i\rho_{i-1}(\sigma_{i-1}-\varepsilon\sigma_i).\label{eq:z}
\end{eqnarray}
Eliminating $\rho_{i-1}(\sigma_{i-1}-\varepsilon\sigma_i)$ in
(\ref{eq:z}) using (\ref{eq:indhyp}) we get
\begin{equation}\label{eq:yplusz}
y+z=-c_i(\sigma_{i+1}-\sigma_{i-1})(\rho_{i-1}+
\varepsilon\rho_i).
\end{equation}
We now simplify $x+y+z$. In (\ref{eq:x}) we eliminate $a_i, k,
\rho$ 
 using (\ref{eq:aiformA}),(\ref{eq:kformA}), (\ref{eq:rho})
 and simplify the result using the fact that $g=h\frac{\varepsilon-1}{1-\sigma}$ from
(\ref{eq:gformA}), (\ref{eq:hformA}); and in (\ref{eq:yplusz}) we
eliminate $c_i(\sigma_{i+1}-\sigma_{i-1})$
 using (\ref{eq:ciformA}).   We routinely obtain
$x+y+z$ is equal to
$h{\frac{\sigma_{i+1}-\sigma\sigma_i}{\sigma_{i-1}-\sigma_i}}$
times $\Delta$ where $\Delta$ is  the left-hand side of
(\ref{eq:indhyp}) minus the right-hand side of
(\ref{eq:indhyp}).  Observe $\Delta=0$ by (\ref{eq:indhyp}) so $x+y+z=0$.  Thus we have 
(\ref{eq:bitimeseq}).  We have now shown by induction that (\ref{eq:indhyp}) holds for
$1 \le i \le D.$  Applying Theorem \ref{thm:main1} we find
$\sigma_0, \sigma_1, \ldots, \sigma_D$ and $\rho_0, \rho_1, \ldots, \rho_D$
   form a tight pair.  
   Comparing (\ref{eq:auxform}), (\ref{eq:rho}) we find $\varepsilon$ is the corresponding
   auxiliary parameter.  Applying Definition \ref{def:tightpseudo} we find
   $\sigma_0, \sigma_1, \ldots, \sigma_D$ is tight and $\varepsilon$ is the corresponding 
   auxiliary parameter.
\end{proof}

\section{Feasible pseudo cosine sequences}

Let $\Gamma$ denote a distance-regular graph with diameter $D\ge 3$ and $a_1 \neq 0$.
Suppose we are given a nontrivial pseudo cosine sequence $\sigma_0, \sigma_1, \ldots, \sigma_D$ which is tight.  In view of Lemma \ref{lem:kbiciform} it is natural to consider the case in which 
$\sigma_{i-1} \not=\sigma_{i+1}$ for $1 \leq i \leq D-1.$  We now investigate this case.

\begin{definition}\label{def:feasible}
{\rm
Let
$\Gamma$ denote a distance-regular graph with diameter
$D\ge 3$ and $a_1\not=0$.
Let  $\sigma_0, \sigma_1, \ldots, \sigma_D$
denote a pseudo cosine sequence. We say this sequence
is \emph{feasible} whenever the following (i),(ii) hold.
\begin{enumerate}
\item
$\sigma_{i-1} \not=\sigma_{i+1}$ for $1 \leq i \leq D-1,$
\item
$\sigma_0, \sigma_1, \ldots, \sigma_D$
is tight.
\end{enumerate}}
\end{definition}

\begin{theorem}
Let
$\Gamma$ denote a distance-regular graph with diameter
$D\ge 3$.
Let $\sigma_0, \sigma_1, \ldots, \sigma_D$ and  $\varepsilon, h$
denote real numbers.
Then the following (i), (ii) are equivalent.
\begin{enumerate}
\item $a_1 \not= 0;$ the sequence $\sigma_0, \sigma_1, \ldots, \sigma_D$ is a
feasible pseudo cosine sequence of $\Gamma$, $\varepsilon$ is the
corresponding auxiliary parameter, and
\begin{equation}\label{eq:hformB}
h = \frac{(1-\sigma)(1-\sigma_2)}{(\sigma^2-\sigma_2)
(1-\varepsilon \sigma)}.
\end{equation}
\item $\sigma_0=1$, $\varepsilon \not=1,$ $\varepsilon \not=-1,$
\begin{eqnarray}
k &=& h\frac{\sigma-\varepsilon}{\sigma - 1}, \label{eq:kformB}
\\
b_i &=& h \frac{(\sigma_{i-1}-\sigma \sigma_i)
(\sigma_{i+1}-\varepsilon \sigma_i)} {(\sigma_{i-1}- \sigma_{i+1})
(\sigma_{i+1}- \sigma_{i})} \qquad (1 \leq i \leq D-1),
\label{eq:biformB}
\\
c_i &=& h \frac{(\sigma_{i+1}-\sigma \sigma_i)
(\sigma_{i-1}-\varepsilon \sigma_i)} {(\sigma_{i+1}- \sigma_{i-1})
(\sigma_{i-1}- \sigma_{i})} \qquad (1 \leq i \leq D-1),
\label{eq:ciformB}
\end{eqnarray}
and the denominators in (\ref{eq:kformB})--(\ref{eq:ciformB}) are
nonzero.
\end{enumerate}
\end{theorem}

\begin{proof}
(i) $\Longrightarrow$ (ii) Observe $\sigma_0=1$ by Lemma \ref{lem:recur2} and since
$\sigma_0, \sigma_1, \ldots, \sigma_D$ is a pseudo cosine sequence.
By Corollary 
\ref{cor:epsgood} we have $\varepsilon  \not\in \lbrace 1, -1 \rbrace$. The remaining assertions follow from Theorem \ref{th:mainth}.\\
(ii) $\Longrightarrow$ (i)  First we show (\ref{eq:hformB}).
Setting $i=1$ in (\ref{eq:ciformB}) and
solving for $h$  using the fact that $c_1=1$ we get
(\ref{eq:hformB}). Observe $h\not= 0;$ otherwise $k=0$.  Next we show that $\sigma_0, \sigma_1, \ldots,
\sigma_D$ is a pseudo cosine sequence.  To do this we apply Lemma
\ref{lem:recur}(i),(ii). Using (\ref{eq:kformB}), (\ref{eq:biformB}),
(\ref{eq:ciformB}) we routinely verify that (\ref{eq:recur})
holds.  By Lemma
\ref{lem:recur} $\sigma_0, \sigma_1, \ldots,
\sigma_D$ is a pseudo cosine sequence.  This sequence is nontrivial;
otherwise $h=0.$ 
 Solving for $a_i$ using (\ref{eq:intnos}) and (\ref{eq:biformB}),
(\ref{eq:ciformB})   we  routinely verify that
(\ref{eq:aiformA}), (\ref{eq:gformA}) hold. We show $a_1
\not=0.$ Setting $i=1$ in (\ref{eq:aiformA}) we
obtain
\begin{equation}
a_1=\frac{(1-\varepsilon)(1+\sigma)(1-\sigma_2)}{(1-\varepsilon\sigma)(\sigma_2-\sigma)}.
\end{equation}
Observe $1-\varepsilon \ne 0$ by assumption; 
$1+\sigma\not=0$ by
(\ref{eq:biformB}) and the fact that $b_1\ne 0$;
and 
$1-\sigma_2 \not = 0$
otherwise $h=0$. It follows $a_1\ne 0$.   
Applying Definition \ref{def:condition}(i) we find 
$\sigma_0, \sigma_1, \ldots,
\sigma_D$ and  $\varepsilon$ together satisfy condition $A$.
Now by Theorem \ref{th:mainth}(i),(ii) we find 
$\sigma_0, \sigma_1, \ldots, \sigma_D$  is tight and that $\varepsilon$ is the 
corresponding auxiliary parameter.  We observe 
 $\sigma_{i-1} \not=\sigma_{i+1} \ (1 \leq i\leq D-1)$ since the denominators  in (\ref{eq:biformB}) are nonzero. Now
 $\sigma_0, \sigma_1, \ldots, \sigma_D$ is feasible by Definition \ref{def:feasible}.
\end{proof}

We finish this section with a comment.

\begin{lemma}
Let $\Gamma$ denote a distance-regular graph with diameter $D\ge
3$ and $a_1 \not=0$.
Let $\sigma_0, \sigma_1, \ldots, \sigma_D$  denote
a feasible pseudo cosine sequence. Let
$\rho_0, \rho_1, \ldots, \rho_D$ denote  a  nontrivial pseudo cosine sequence such that
$\sigma_0, \sigma_1, \ldots, \sigma_D$ and $\rho_0, \rho_1, \ldots, \rho_D$
   form a tight pair.  Let $\varepsilon$ denote the corresponding auxiliary parameter.
Then 
\begin{equation}\label{eq:rhoi} \rho_i = {\prod_{j=1}^i}
\frac{\sigma_{j-1}-\varepsilon\sigma_j}{\sigma_j-\varepsilon\sigma_{j-1}}\qquad
(0 \le i \le D) \end{equation}  and the denominators in
(\ref{eq:rhoi}) are all nonzero. 
\end{lemma}

\begin{proof}
We first show the denominators in (\ref{eq:rhoi})
are all nonzero.
By Lemma \ref{lem:zerorhoi} and since  
$\sigma_{i-1}-\sigma_{i+1}\not=0$ for $1 \leq i \leq D-1$ we find
 $\sigma_i-\varepsilon\sigma_{i-1}\not=0$
for $2\le i \le D$. Observe $\sigma\not=\varepsilon$ in view of Lemma
\ref{lem:rhorho2}.  Using this and the fact $\sigma_0=1$ we get
$\sigma_i-\varepsilon\sigma_{i-1}\not=0$ for $1 \le i \le D$.
We have now shown the denominators in (\ref{eq:rhoi}) are all nonzero.
We now verify (\ref{eq:rhoi}).  Observe (\ref{eq:rhoi}) holds for $i=0$
since $\rho_0=1$.  Line (\ref{eq:rhoi}) holds for $i=1$ by 
(\ref{eq:auxform}).  Line (\ref{eq:rhoi}) holds for $2 \leq i \leq D$ by
(\ref{eq:rhoplus}) and a routine induction. 
\end{proof}

\noindent Arlene A. Pascasio \\
Department of Mathematics \\
De La Salle University - Manila \\
2401 Taft Avenue\\
Malate Manila 1004\\
Philippines\\
Email: pascasioa@dlsu.edu.ph\\

\noindent Paul Terwilliger \\
 Department of Mathematics \\ 
University of Wisconsin - Madison \\
Van Vleck Hall \\
480 Lincoln Drive \\
Madison WI, USA  53706-1388 \\
Email: terwilli@math.wisc.edu


\begin{thebibliography}{9}

\bibitem{bi} E. Bannai and T. Ito, {\it Algebraic Combinatorics I:
Association Schemes}, Benjamin-Cummings Lecture Note Ser. 58,
Benjamin-Cummings, Menlo Park, CA 1984.

\bibitem{bcn} A.E. Brouwer, A.M.
Cohen and A. Neumaier, {\it Distance-Regular Graphs}, Springer, New York,
1989.
 
\bibitem{go2} J. T. Go and P. Terwilliger,
Tight distance-regular graphs and the subconstituent algebra,
European J. Combin., 23 (2002), 793--816.

\bibitem{jk5}  A. Juri\v si\'c, AT4 family and 
2-homogeneous graphs, Preprint.

\bibitem{jk}  A. Juri\v si\'c and J. Koolen,  A local approach
to 1-homogeneous graphs, Des. Codes Cryptography, 
21 (2000), 127--147.

\bibitem{jk0}  A. Juri\v si\'c and J. Koolen,  Nonexistence 
of some antipodal distance-regular graphs of diameter four,
European J. Combin., 21 (2000), 1039--1046.

\bibitem{jk2}  A. Juri\v si\'c and J. Koolen,  1-homogeneous
graphs with cocktail party $\mu$-graphs, J. Algebraic Combin., To appear. 

\bibitem{jk3}  A. Juri\v si\'c and J. Koolen,  Krein parameters
and antipodal distance-regular graphs with diameter 3 and 4,
Discrete Math., 244 (2002), 181--202. 


\bibitem{jkt} A. Juri\v si\'c, J. Koolen and P. Terwilliger,
Tight distance-regular graphs, 
J. Algebraic Combin., 12 (2000), 163--197. 

\bibitem{jurter} A. Juri\v si\'c and P. Terwilliger,
Pseudo tight distance-regular graphs, Preprint.


\bibitem{maclean1} M. MacLean,  An inequality involving two eigenvalues of a bipartite
distance-regular graph, Discrete Math., 225 (2000), 193--216.

\bibitem{maclean2} M. MacLean, Taut distance-regular graphs
of odd diameter, J. Algebraic Combin., Submitted.

\bibitem{maclean3} M. MacLean,  Taut distance-regular graphs
of even diameter, J. Combin. Theory Ser. B, Submitted. 


\bibitem{aap1} A. A. Pascasio, Tight graphs and their primitive
idempotents, J. Algebraic Combin., 10 (1999), 47-59.

\bibitem{aap2} A. A. Pascasio,  Tight distance-regular
graphs and Q-polynomial property, Graphs and Combin., 
    17 (2001), 149--169.

\bibitem{aap3} A. A. Pascasio, An inequality on the cosines
of a tight distance-regular graph, Linear Algebra Appl., 
    325 (2001), 147--159.

\bibitem{aap4} A. A. Pascasio,
An inequality in character algebras, Discrete Math., 264 (2003), 201--209.



\bibitem{tomiyama} M. Tomiyama, On the primitive idempotents
of distance-regular graphs, Discrete Math., 240 (2001), 281--294.

\bibitem{terweng} C.W. Weng and P. Terwilliger,
Distance-regular graphs, pseudo primitive idempotents,
and the Terwilliger algebra, Seidel 80 Conference, To appear.



\end{thebibliography}
\end{document}